\documentclass{amsart}
\usepackage{amssymb,amsmath, amsthm,latexsym}
\usepackage{graphics}
\usepackage{amscd}
\usepackage{graphics}
\usepackage{here}
\newcommand{\cal}[1]{\mathcal{#1}}
\theoremstyle{plain}
\newtheorem*{theo}{Theorem}

\newtheorem*{lem}{Lemma}
\newtheorem*{prop}{Proposition}
\newtheorem{lemma}{Lemma}[section]
\newtheorem{theorem}[lemma]{Theorem}
\newtheorem{proposition}[lemma]{Proposition}
\newtheorem{corollary}[lemma]{Corollary}
\theoremstyle{definition}
\newtheorem*{definition}{Definition}
\parskip=\bigskipamount

\let\egthree=\phi
\let\phi=\varphi
\let\varphi=\egthree

\begin{document}
\title{Invariant Radon measures on measured lamination space}
\author{Ursula Hamenst\"adt}
\thanks
{
Research
partially supported by DFG Sonderforschungsbereich 611}

\begin{abstract}
Let $S$ be an oriented surface of
genus $g\geq 0$ with $m\geq 0$ punctures
and $3g-3+m\geq 2$. We classify
all Radon measures on the space
of measured geodesic laminations which are
invariant under the action of the mapping class
group of $S$.
\end{abstract}

\maketitle

\section{Introduction}

Let $S$ be an oriented surface of finite type, i.e. $S$ is a
closed surface of genus $g\geq 0$ from which $m\geq 0$
points, so-called \emph{punctures},
have been deleted. We assume that $3g-3+m\geq 1$,
i.e. that $S$ is not a sphere with at most $3$
punctures or a torus without puncture.
In particular, the Euler characteristic of $S$ is negative.
Then the  \emph{Teichm\"uller space} ${\cal T}(S)$
of $S$ is the quotient of the space of all complete
hyperbolic metrics of finite volume on $S$ under the action of the
group of diffeomorphisms of $S$ which are isotopic
to the identity.
The \emph{mapping class group} ${\cal M}(S)$
of all isotopy classes
of orientation preserving diffeomorphisms
of $S$ acts properly discontinuously
on ${\cal T}(S)$ with quotient the
\emph{moduli space} ${\rm Mod}(S)$.

A \emph{geodesic lamination} for a fixed choice of a complete
hyperbolic metric of finite volume on $S$
is a \emph{compact} subset of $S$ foliated
into \emph{simple} geodesics. A \emph{measured geodesic
lamination} is a geodesic lamination together with
a transverse translation invariant measure.
The space ${\cal M\cal L}$ of all measured
geodesic laminations on $S$, equipped
with the weak$^*$-topology,
is homeomorphic to
$S^{6g-7+2m}\times (0,\infty)$ where $S^{6g-7+2m}$ is the
$6g-7+2m$-dimensional sphere.
The mapping class group ${\cal M}(S)$
naturally acts on ${\cal M\cal L}$ as a group of
homeomorphisms preserving a Radon measure
in the \emph{Lebesgue measure class}.
Up to scale, this measure is induced
by a natural symplectic structure on ${\cal M\cal L}$
(see \cite{PH92} for this observation of
Thurston), and it is ergodic under the action of
${\cal M}(S)$. Moreover, the measure is
\emph{non-wandering}. By this we mean that
${\cal M\cal L}$ does
not admit an ${\cal M}(S)$-invariant
countable Borel partition into sets
of positive measure.

If $S$ is a once punctured torus or a
sphere with 4 punctures, then
the Teichm\"uller space of $S$ has a natural
identification with the upper half-plane
${\bf H}^2=\{z\in \mathbb{C}\mid {\rm Im}(z)>0\}$.
Up to passing to the quotient by the
hyperelliptic involution if $S$ is the once punctured torus,
the mapping class group ${\cal M}(S)$ is just
the group $PSL(2,\mathbb{Z})$ acting on ${\bf H}^2$ by linear fractional
transformations. The
action of ${\cal M}(S)$ on
measured lamination space can in this case be identified
with the quotient of the standard
\emph{linear} action of  $SL(2,\mathbb{Z})$ on
$\mathbb{R}^2$ under the reflection at the origin (see
the book \cite{BM00} of Bekka and Mayer for
more and for references and compare the survey \cite{H07c}).

Extending earlier work of
Furstenberg, in 1978 Dani
completely classified all $SL(2,\mathbb{Z})$-invariant
Radon measures on $\mathbb{R}^2$. He showed
that if such a measure $\eta$ is ergodic under the
action of $SL(2,\mathbb{Z})$ then either it
is non-wandering and coincides with the usual
Lebesgue measure $\lambda$ up to scale, or it is
\emph{rational}, which means that it is
supported on a single $SL(2,\mathbb{Z})$-orbit
of points whose coordinates are dependent over $\mathbb{Q}$.

If the surface $S$ is \emph{non-exceptional},
i.e. if $3g-3+m\geq 2$, then the
${\cal M}(S)$-invariant Radon measures on ${\cal M\cal L}$ which
naturally
correspond to the rational measures for exceptional
surfaces are defined as follows.
A \emph{weighted geodesic multi-curve} on $S$ is
a measured geodesic lamination whose \emph{support}
is a union of simple closed geodesics.
The orbit of a weighted geodesic multi-curve under
the action of ${\cal M}(S)$ is a discrete subset of
${\cal M\cal L}$ (see Section 5 for this easy and well known
fact) and hence it supports a ray
of invariant purely atomic Radon measures which we
call \emph{rational}. 
This definition coincides with the one
given above for a once punctured torus or
a forth punctured sphere.

For a non-exceptional surface $S$,
there are additional ${\cal M}(S)$-invariant
Radon measures on ${\cal M\cal L}$.
Namely, a \emph{proper connected bordered subsurface} $S_0$
of $S$ is a connected component of
the space which we obtain from $S$ by cutting $S$ open along
a collection of disjoint simple closed geodesics.
Then $S_0$ is a connected
surface with non-empty geodesic boundary and of
negative Euler characteristic.
If two boundary components of $S_0$
correspond to the same closed geodesic $\gamma$ in $S$ then
we require that $S-S_0$ contains a connected component which
is an annulus with core curve $\gamma$.
Let ${\cal M\cal L}(S_0)\subset {\cal M\cal L}$ be the space
of all measured geodesic laminations on $S$ which
are contained in the interior of $S_0$.
The space ${\cal M\cal L}(S_0)$ can naturally be
identified with the space of measured geodesic
laminations on the
surface $\hat S_0$ of finite type which we obtain
from $S_0$ by collapsing each boundary circle to a puncture.
The stabilizer in ${\cal M}(S)$ of the subsurface $S_0$
is the direct product of the group of all elements which can be
represented by diffeomorphisms leaving
$S_0$ pointwise fixed and a group which is naturally isomorphic to
a subgroup of finite index of the mapping class
group ${\cal M}(\hat S_0)$ of
$\hat S_0$.

Let $c$ be a
weighted geodesic multi-curve on $S$ which
is disjoint from the interior of $S_0$. Then for every
$\zeta\in {\cal M\cal L}(S_0)$ the union
$c\cup \zeta$ is a measured geodesic lamination on $S$
which we denote by $c\times \zeta$.
Let $\mu(S_0)$ be an
${\cal M}(\hat S_0)$-invariant Radon measure on ${\cal M\cal L}(S_0)$
which is contained in the Lebesgue measure class.
The measure $\mu(S_0)$ can
be viewed as a Radon measure on ${\cal M\cal L}$ which
gives full measure to the laminations of the form
$c\times \zeta$ $(\zeta\in {\cal M\cal L}(S_0))$ and which
is invariant and ergodic under the stabilizer of
$c\cup S_0$ in ${\cal M}(S)$.
The translates of this measure under the
action of ${\cal M}(S)$ define an ${\cal M}(S)$-invariant
ergodic wandering measure on ${\cal M\cal L}$
which we call a \emph{standard subsurface measure}.
We observe in Section 5 that if
the weighted geodesic multi-curve
$c$ contains the boundary of $S_0$ then the standard
subsurface measure defined by $\mu(S_0)$ and $c$
is a Radon measure on ${\cal M \cal L}$.

The goal of this note is to show that
every ${\cal M}(S)$-invariant ergodic Radon measure
on ${\cal M\cal L}$ is of the form described above.

\begin{theo}
\begin{enumerate}
\item
An invariant ergodic non-wandering Radon measure for
the action of ${\cal M}(S)$ on ${\cal M\cal L}$
coincides with the Lebesgue measure up to scale.
\item An invariant ergodic wandering
Radon measure for the action of ${\cal M}(S)$
on ${\cal M\cal L}$ is either rational or a standard
subsurface measure.
\end{enumerate}
\end{theo}

The organization of the paper is as follows.
In Section 2 we discuss some properties of
geodesic laminations, quadratic differentials
and the curve graph needed in the sequel.
In Section 3 we introduce conformal densities for the
mapping class group. These
conformal densities are families of finite Borel measures
on the projectivization ${\cal P\cal M\cal L}$
of ${\cal M\cal L}$, parametrized
by the points in Teichm\"uller space.
They are defined in analogy to
the conformal densities for discrete subgroups of
the isometry group of a hyperbolic space. Up to scale,
there is a unique conformal density in the Lebesgue measure
class. We show that this is the only conformal density which
gives full measure to the ${\cal M}(S)$-invariant subset
of ${\cal P\cal M\cal L}$ of all
projective measured geodesic
laminations whose support is minimal and fills up $S$.

Every conformal density gives rise to an ${\cal M}(S)$-invariant
Radon measure on ${\cal M\cal L}$. The investigation of
invariant Radon measure which are not of this form relies
on the structural results of Sarig \cite{S04}.
To apply his results we use train tracks to
construct partitions of measured lamination space
which have properties similar to Markov partitions.
Section 4 summarizes those
facts about train tracks which are needed for this purpose.
The proof of the theorem is completed in Section 5.
In the appendix we present
a result of Minsky and Weiss
\cite{MW02} in the form needed for the proof of our
theorem.

After this paper was posted on the arXiv I received
the preprint \cite{LM07} of Lindenstrauss and
Mirkzakhani which contains another proof of the
above theorem.

\section{Quadratic differentials and the curve graph}

In this introductory section we summarize some
properties of (measured) geodesic laminations,
quadratic differentials and the curve graph which are needed later on.
We also introduce some notations which will be used throughout the paper.

In the sequel we always denote by $S$ an
oriented surface of
genus $g\geq 0$ with $m\geq 0$ punctures and where $3g-3+m\geq 1$.

\subsection{Geodesic laminations}

A \emph{geodesic lamination} for a complete
hyperbolic structure of finite volume on the surface $S$ is
a \emph{compact} subset of $S$ which is foliated into simple
geodesics.
A geodesic lamination $\lambda$ is called \emph{minimal}
if each of its half-leaves is dense in $\lambda$. Thus a simple
closed geodesic is a minimal geodesic lamination. A minimal
geodesic lamination with more than one leaf has uncountably
many leaves and is called \emph{minimal arational}.
Every geodesic lamination $\lambda$ consists of a disjoint union of
finitely many minimal components and a finite number of isolated
leaves. Each of the isolated leaves of $\lambda$ either is an
isolated closed geodesic and hence a minimal component, or it
\emph{spirals} about one or two minimal components
\cite{CEG87}.
A geodesic lamination is \emph{maximal}
if its complementary regions are all ideal triangles
or once punctured monogons. A geodesic lamination
\emph{fills up $S$} if its complementary
regions are all topological discs or once
punctured topological discs.

A \emph{measured geodesic lamination} is a geodesic lamination
$\lambda$ together with a translation invariant transverse
measure. Such a measure assigns a positive weight to each compact
arc in $S$ which intersects $\lambda$ nontrivially and
transversely and whose endpoints are contained
in the complementary regions of
$\lambda$. The geodesic lamination $\lambda$ is called the
\emph{support} of the measured geodesic lamination; it consists of
a disjoint union of minimal components \cite{CEG87}.
The space ${\cal M\cal L}$
of all measured geodesic laminations on $S$ equipped with the
weak$^*$-topology is homeomorphic to $S^{6g-7+2m}\times
(0,\infty)\sim \mathbb{R}^{6g-6+2m}-\{0\}$.
Its projectivization is the space ${\cal P\cal M\cal L}$
of all \emph{projective measured geodesic laminations}.
The measured geodesic lamination $\mu\in {\cal
M\cal L}$ \emph{fills up $S$} if its support fills up $S$.
The projectivization of a measured geodesic lamination
which fills up $S$ is also said to fill up $S$.
A measured geodesic lamination is called
\emph{uniquely ergodic} if its support admits a single
transverse measure up to scale.
There is a continuous symmetric pairing $i:{\cal M\cal L}\times
{\cal M\cal L}\to [0,\infty)$, the so-called \emph{intersection form},
which extends the geometric intersection number between two
simple closed curves.

\subsection{Quadratic differentials}

The smooth
fibre bundle ${\cal Q}^1(S)$
of all \emph{holomorphic
quadratic differentials of area
one} over the Teichm\"uller space ${\cal T}(S)$ of the surface $S$
can naturally be viewed as the unit cotangent
bundle of ${\cal T}(S)$ for the \emph{Teichm\"uller metric}.
The \emph{Teichm\"uller geodesic flow} $\Phi^t$ on
${\cal Q}^1(S)$ commutes
with the action of the mapping class group
${\cal M}(S)$ of all isotopy classes of
orientation preserving diffeomorphisms of $S$.
Thus this flow descends to a flow
on the quotient ${\cal Q}(S)=
{\cal Q}^1(S)/{\cal M}(S)$, again denoted
by $\Phi^t$.

A measured geodesic lamination can be viewed
as an equivalence class of \emph{measured foliations} on $S$
\cite{L83}.
Therefore every holomorphic quadratic differential
$q$ on $S$ defines a pair $(\mu,\nu)\in {\cal M\cal L}\times
{\cal M\cal L}$ where the \emph{horizontal measured
geodesic lamination} $\mu$
corresponds to the \emph{horizontal
measured foliation} of $q$ which is expanded under the
Teichm\"uller flow, and where the
\emph{vertical measured geodesic lamination} $\nu$
corresponds to the vertical measured foliation of $q$
which is contracted under the Teichm\"uller flow.
The area of the quadratic differential is just the intersection
$i(\mu,\nu)$.

For a quadratic differential $q\in {\cal Q}^1(S)$ define the
\emph{strong unstable manifold} $W^{su}(q)$
to be the set of all quadratic differentials $z\in {\cal Q}^1(S)$
whose vertical measured geodesic lamination equals
the vertical measured geodesic lamination
for $q$.
Similarly, define the \emph{strong stable
manifold} $W^{ss}(q)$ to be the set of all quadratic
differentials $z\in {\cal Q}^1(S)$
whose horizontal measured geodesic
lamination coincides with the horizontal measured
geodesic lamination of $q$.
The \emph{stable manifold}
$W^s(q)=\cup_{t\in \mathbb{R}}\Phi^tW^{ss}(q)$  and
the \emph{unstable manifold}
$W^u(q)=\cup_{t\in \mathbb{R}}\Phi^t W^{su}(q)$ are
submanifolds
of ${\cal Q}^1(S)$ which project homeomorphically
onto ${\cal T}(S)$ \cite{HM79}.

The sets $W^s(q)$ (or $W^{ss}(q),W^{su}(q),W^u(q)$)
$(q\in {\cal Q}^1(S))$
define a foliation of
${\cal Q}^1(S)$ which is invariant under the mapping class group
and hence projects to a singular
foliation on ${\cal Q}(S)$
which we call the \emph{stable foliation}
(or the \emph{strong stable, strong unstable, unstable foliation}).
There is a distinguished family of Lebesgue measures
$\lambda^s$ on the leaves of the stable
foliation which are conditional measures of a
$\Phi^t$-invariant Borel probability measure
$\lambda$ on ${\cal Q}(S)$ in the Lebesgue measure
class.
The measure $\lambda$ is ergodic and mixing
under the Teichm\"uller geodesic flow
(see \cite{M82a} and also \cite{V82,V86}).

Every area one quadratic differential $q\in {\cal Q}^1(S)$
defines a singular euclidean metric on $S$ of area one
together with two orthogonal foliations
by straight lines, with singularities
of cone angle $k\pi$ for some $k\geq 3$.
This metric is given
by a distinguished
family of isometric charts $\phi:U\subset S\to
\phi(U)\subset \mathbb{C}$ on the complement of
the zeros (or poles at the punctures) of $q$ which
map the distinguished foliations to the
foliation of $\mathbb{C}$ into the \emph{horizontal}
lines parallel to the
real axis and into the \emph{vertical lines} parallel to
the imaginary axis.

Let $\hat S$ be the compactification of $S$ which
we obtain by adding a point at each puncture.
A \emph{saddle connection} for $q$ is a path
$\delta:(0,1)\to S$ which does  not contain any
singular point, whose image under a distinguished chart
is an euclidean straight line and which extends continuously
to a path $\bar\delta:[0,1]\to \hat S$ mapping the
endpoints to singular points or punctures.
A saddle connection
is \emph{horizontal} if it is mapped by a distinguished
chart to a horizontal line segment.

The group $SL(2,\mathbb{R})$ acts on
the bundle ${\cal Q}^1(S)$ by replacing for
$q\in {\cal Q}^1(S)$ and $M\in SL(2,\mathbb{R})$
each isometric chart $\phi$
for $q$ by $M\circ \phi$ where $M$ acts linearly
on $\mathbb{R}^2=\mathbb{C}$. This preserves the compatibility
condition for charts. The Teichm\"uller geodesic flow $\Phi^t$
then is the action of the diagonal group
\begin{equation} \Bigr( \begin{matrix}e^{t} & 0\\
0 & e^{-t}\\
\end{matrix} \Bigl) \quad  (t\in \mathbb{R}) .
\end{equation}
The so-called \emph{horocycle flow} $h_t$ is given by the action of the
unipotent subgroup
\begin{equation} \Bigr( \begin{matrix}
1 & 0\\
t & 1 \\
\end{matrix} \Bigl)  \quad (t\in \mathbb{R}) .
\end{equation}

\subsection{The curve graph}

The \emph{curve graph} ${\cal C}(S)$ of the surface $S$
is a metric graph whose vertices are the
free homotopy classes of \emph{essential}
simple closed
curves on $S$, i.e. curves which are neither
contractible nor freely homotopic into a puncture
of $S$. In the sequel we often do not distinguish between an
essential simple closed curve and its free
homotopy class whenever no confusion is possible.
Two such curves
are connected in ${\cal C}(S)$
by an edge of length one
if and only if they can be realized
disjointly. If the surface $S$ is non-exceptional,
i.e. if $3g-3+m\geq 2$, then
the curve graph ${\cal C}(S)$ is connected.
Any two elements $c,d\in {\cal C}(S)$
of distance at least 3
\emph{jointly fill up $S$}, i.e.
they decompose $S$ into topological discs
and once punctured topological discs.
The mapping class group naturally acts
on ${\cal C}(S)$ as a group of simplicial
isometries.

By Bers' theorem, there is a number $\chi_0>0$
such that for every complete hyperbolic metric
$h$ on $S$ of finite volume there is a \emph{pants decomposition}
of $S$ consisting of $3g-3+m$ pairwise disjoint
simple closed geodesics
of length at most $\chi_0$.
On the other hand, the number of essential
simple closed curves $\alpha$
on $S$ whose hyperbolic length
$\ell_h(\alpha)$ (i.e. the length of a geodesic representative
of its free homotopy class)
does not exceed $2\chi_0$ is bounded from above by a
universal constant not depending on $h$, and the diameter of the
subset of ${\cal C}(S)$ containing these curves is uniformly
bounded as well.

Define a map \begin{equation}
\Upsilon_{\cal T}:{\cal T}(S)\to {\cal C}(S)
\end{equation} by
associating to a complete hyperbolic metric $h$ on $S$
of finite volume a curve $\Upsilon_{\cal T}(h)$ whose $h$-length
is at most $\chi_0$. If $\Upsilon^\prime$ is any other choice of
such a map then $d(\Upsilon_{\cal T}(h),\Upsilon^\prime_{\cal T}(h)\leq
{\rm const}$.
By the discussion in \cite{H06a} there
is a constant $L>1$ such that
\begin{equation}\label{Lipschitz}
d(\Upsilon_{\cal T}(g),\Upsilon_{\cal T}(h))
\leq L d(g,h) +L\quad \text{for all} \quad g,h\in {\cal T}(S)
\end{equation}
where by abuse, we use the same symbol $d$ to denote the
distance on ${\cal T}(S)$ defined by the Teichm\"uller metric
and the distance on the curve graph ${\cal C}(S)$.

For a quadratic differential $q\in {\cal Q}^1(S)$ define
the \emph{$q$-length} of an essential
closed curve $\alpha$ on $S$ to be the
minimal length of a representative of the
free homotopy class of $\alpha$ with
respect to the singular euclidean metric
defined by $q$. We have.

\begin{lemma}\label{simple1} For every $\chi>0$ there
is a constant $a(\chi)>0$ with the
following property.
For any quadratic
differential $q\in {\cal Q}^1(S)$
the diameter in ${\cal C}(S)$ of the set of
all simple closed curves on $S$ of $q$-length at
most $\chi$ does not exceed $a(\chi)$.
\end{lemma}

\begin{proof} By Lemma 5.1 of \cite{MM99}
(see also Lemma 5.1 of \cite{Bw06} for an alternative
proof)
there is a number $b>0$ and for every
singular euclidean metric on $S$ defined
by a quadratic differential $q$ of area one
there is an embedded annulus $A\subset S$ of
width at least $b$. If we denote by $\gamma$
the core curve of $A$, then
the $q$-length of every simple closed curve
$c$ on $S$ is at least $i(c,\gamma)b$.
As a consequence, for every $\chi >0$
the intersection number
with $\gamma$ of every simple closed
curve $c$ on $S$ whose $q$-length is at most
$\chi$ is bounded from above by $\chi/b$.
This implies that the
set of such curves is of uniformly bounded
diameter in ${\cal C}(S)$ (see \cite{MM99,Bw06}).
\end{proof}

By possibly enlarging the constant
$\chi_0>0$ as above we may assume that
for every $q\in {\cal Q}^1(S)$
there is an essential simple closed
curve on $S$ of $q$-length
at most $\chi_0$ (see the
proof of Lemma \ref{simple2} for a justification
of this well known fact). Thus we
can define a map
\begin{equation}\label{upsilonq}
\Upsilon_{\cal Q}: {\cal Q}^1(S)\to
{\cal C}(S)
\end{equation} by associating to
a quadratic differential $q$ a simple
closed curve $\Upsilon_{\cal Q}(q)$ whose
$q$-length is at most $\chi_0$.
By Lemma \ref{simple1}, if
$\Upsilon^\prime_{\cal Q}$ is any other choice of
such a map then we have $d(\Upsilon_{\cal Q}(q),
\Upsilon_{\cal Q}^\prime(q))\leq
a(\chi_0)$ for all $q\in {\cal Q}^1(S)$ where
$d$ is the distance function on ${\cal C}(S)$.

Recall the definition of the map
$\Upsilon_{\cal T}:{\cal T}(S)\to {\cal C}(S)$
which associates to a complete hyperbolic metric $h$
on $S$ of finite volume a simple closed curve
$c$ of $h$-length at most $\chi_0$.
Let $P:{\cal Q}^1(S)\to {\cal T}(S)$ be
the canonical projection.
The following simple
observation is related to recent
work of Rafi \cite{R05}. We include its short proof for
completeness.

\begin{lemma}\label{simple2} There is a constant
$\chi_1>0$ such that $d(\Upsilon_{\cal Q}(q),
\Upsilon_{\cal T}(Pq))\leq \chi_1$ for all
$q\in {\cal Q}^1(S)$.
\end{lemma}

\begin{proof} By Lemma \ref{simple1},
it is enough to show that
for every $q\in {\cal Q}^1(S)$ and every simple
closed curve $\alpha$ on $S$ whose
length with respect to the hyperbolic metric $Pq$
is bounded from above by $\chi_0$,
the $q$-length of $\alpha$ is uniformly
bounded.

For this observe that by the collar lemma of
hyperbolic geometry,
a simple closed geodesic $\alpha$
on a hyperbolic surface
whose length is bounded from above by
$\chi_0$ is the core curve
of an embedded annulus $A$ whose \emph{modulus} is bounded
from below by a universal constant $b>0$.
Then the \emph{extremal length} of the
core curve of $A$ is bounded from above
by a universal constant $c>0$. Now the area of
$q$ equals one and therefore the $q$-length
of the core curve $\alpha$ of $A$ does not
exceed $\sqrt{c}$ by the definition of extremal
length (see e.g. \cite{M96}). This shows the lemma.
\end{proof}

Choose a smooth function $\sigma:[0,\infty)\to [0,1]$ with
$\sigma[0,\chi_0]\equiv 1$ and
$\sigma[2\chi_0,\infty)\equiv 0$ where as before,
$\chi_0$ is a Bers constant for $S$.
For every $h\in {\cal T}(S)$ we obtain a
finite Borel measure $\mu_h$ on
the curve graph ${\cal C}(S)$ by defining
\begin{equation}
\mu_h=\sum_{c\in {\cal C}(S)} \sigma(\ell_h(c))\delta_c
\end{equation}
where $\delta_c$ denotes the Dirac mass at $c$. The total mass
of $\mu_h$ is bounded from above and below by a universal positive
constant, and the diameter of the support of $\mu_h$ in ${\cal
C}(S)$ is uniformly bounded as well. Moreover, the measures
$\mu_h$ depend continuously on $h\in {\cal T}(S)$ in the
weak$^*$-topology. This means that for every bounded function
$f:{\cal C}(S)\to \mathbb{R}$ the function $h\to \int f d\mu_h$ is
continuous.

The curve graph ${\cal C}(S)$ is a hyperbolic geodesic
metric space \cite{MM99} and hence it admits
a \emph{Gromov boundary} $\partial {\cal C}(S)$.
For every $c\in {\cal C}(S)$ there is a complete distance
function $\delta_c$ on the Gromov boundary $\partial{\cal
C}(S)$ of ${\cal C}(S)$
of uniformly bounded diameter and there is
a number $\beta >0$ such that $\delta_c\leq e^{\beta
d(c,a)}\delta_a$ for all $c,a\in {\cal C}(S)$. 
The distances $\delta_c$ are equivariant with respect to
the action of ${\cal M}(S)$ on ${\cal C}(S)$ and
on $\partial {\cal C}(S)$.

For $h\in {\cal
T}(S)$ define a distance $\delta_h$ on $\partial {\cal C}(S)$ by
\begin{equation}\label{distance}
\delta_h(\xi,\zeta)=\int \delta_c(\xi,\zeta)
d\mu_h(c). \end{equation}
Clearly the metrics
$\delta_h$ are equivariant with respect to the
action of ${\cal M}(S)$ on ${\cal T}(S)$ and
$\partial{\cal C}(S)$. Moreover, there is a constant
$\kappa >0$ such that
\begin{equation}\label{delta}
\delta_h\leq e^{\kappa d(h,z)}\delta_z
\end{equation}
for all $h,z\in {\cal T}(S)$ (here as before, $d$ denotes
the Teichm\"uller metric). Namely, the function $\sigma$ is
smooth, with uniformly bounded differential. Moreover,
for every simple closed curve $c\in {\cal C}(S)$, the
function $h\to \log \ell_h(c)$ on ${\cal T}(S)$ is
smooth, with unformly bounded differential with respect
to the norm induced by the Teichm\"uller metric (see
\cite{IT99}). Since $\sigma$ is supported in $[0,2\chi_0]$,
this implies that for each $c\in {\cal C}(S)$
the function $h\to \sigma(\ell_h(c))$ on ${\cal T}(S)$ is
smooth, with uniformly bounded differential. As a consequence,
for all $\xi\not=\eta\in \partial {\cal C}(S)$ the function
$h\to \delta_h(\xi,\eta)$ is smooth, and the differential
of its logarithm is uniformly bounded with respect
to the Teichm\"uller norm, independent of $\xi,\eta$.
From this
and the definitions, the estimate (\ref{delta}) above
is immediate. Via enlarging the constant $\kappa$
we may also assume that
\begin{equation}\label{deltacomparison}
\kappa^{-1}\delta_h\leq
\delta_{\Upsilon(h)}\leq \kappa\delta_h
\end{equation} for every
$h\in {\cal T}(S)$.

\section{Conformal densities}

In this section we study
\emph{conformal densities} on the space
${\cal P\cal M\cal L}$
of projective measured geodesic laminations on $S$.
Recall that ${\cal P\cal M\cal L}$ equipped with
the weak topology is homeomorphic to the
sphere $S^{6g-7+2m}$, and the
mapping class group ${\cal M}(S)$
naturally acts on ${\cal P\cal M\cal L}$
as a group of homeomorphisms. By the
Hubbard Masur theorem \cite{HM79}, for every
$x\in {\cal T}(S)$
and every $\lambda\in {\cal P\cal M\cal L}$ there
is a unique holomorphic quadratic
differential $q(x,\lambda)\in {\cal Q}^1(S)_x$
of area one on $x$ whose horizontal measured geodesic
lamination $q_h(x,\lambda)$
is contained in the class $\lambda$.
For all $x,y\in {\cal T}(S)$
there is a number $\Psi(x,y,\lambda)\in \mathbb{R}$
such that $q_h(y,\lambda)=e^{\Psi(x,y,\lambda)}q_h(x,\lambda)$.
The function $\Psi:{\cal T}(S)\times {\cal T}(S)\times
{\cal P\cal M\cal L}\to \mathbb{R}$ is continuous, moreover
it satisfies the \emph{cocycle identity}
\begin{equation}\label{cocycle}
\Psi(x,y,\lambda)+\Psi(y,z,\lambda)=\Psi(x,z,\lambda)
\end{equation}
for all $x,y,z\in {\cal T}(S)$ and all $\lambda\in
{\cal P\cal M\cal L}$.

\begin{definition}\label{condense}
A \emph{conformal density of dimension $\alpha\geq 0$} on
${\cal P\cal M\cal L}$ is an ${\cal M}(S)$-equivariant
family $\{\nu^y\}$ $(y\in {\cal T}(S))$ of finite Borel measures
on ${\cal P\cal M\cal L}$ which are absolutely
continuous and satisfy $d\nu^{z}/d\nu^y=e^{\alpha\Psi(y,z,\cdot)}$
almost everywhere.
\end{definition}

The conformal density $\{\nu^y\}$ is \emph{ergodic}
if the ${\cal M}(S)$-invariant measure class it defines
on ${\cal P\cal M\cal L}$ is ergodic.
There is an ergodic
conformal density $\{\lambda^x\}$
of dimension $\alpha=6g-6+2m$ in the Lebesgue measure
class induced by
the symplectic form on the space
${\cal M\cal L}$ of all measured geodesic laminations on $S$, see \cite{M82a}.
Note that since the action of ${\cal M}(S)$ on ${\cal P\cal M\cal L}$ is
minimal, the measure class of a conformal density
is always of full support.

In Section 5 we will see that every conformal density gives
rise to an ${\cal M}(S)$-invariant Radon measure on
${\cal M\cal L}$ (see \cite{S04}), so the classification of
conformal densities is essential for the classification of
invariant Radon measures on ${\cal M\cal L}$.

Following \cite{Su79}, we construct from a
conformal density
$\{\nu^x\}$ of dimension $\alpha$ an
${\cal M}(S)$-invariant family $\{\nu^{su}\}$ of
locally finite Borel measures on strong unstable
manifolds $W^{su}(q)$ $(q\in {\cal Q}^1(S))$
which transform under the Teichm\"uller geodesic flow $\Phi^t$
via $\nu^{su}\circ \Phi^t=
e^{\alpha t}\nu^{su}$. For this let
\begin{equation}\label{pi}
\pi:{\cal Q}^1(S)\to {\cal P\cal M\cal L}
\end{equation}  be
the natural projection which maps a quadratic
differential
$q\in {\cal Q}^1(S)$ to its horizontal projective
measured geodesic lamination.
Let $P:{\cal Q}^1(S)\to {\cal T}(S)$ be the
canonical projection.
The restriction
of the projection $\pi:{\cal Q}^1(S)\to {\cal P\cal M\cal L}$
to the strong unstable manifold
$W^{su}(q)$ is a homeomorphism
onto its image \cite{HM79} and hence the measure
$\nu^{Pq}$ on ${\cal P\cal M\cal L}$
induces a Borel measure
$\tilde \nu^{su}$ on $W^{su}(q)$.
The measure $\nu^{su}$ on $W^{su}(q)$ defined by
$d\nu^{su}(u)=e^{\alpha\Psi(Pq,Pu,\pi(u))}d\tilde \nu^{su}(u)$ is
locally finite and does not depend on the
choice of $q$. The measures $\nu^{su}$ on
strong unstable manifolds transform under the
Teichm\"uller flow as required.

The \emph{flip} ${\cal F}:q\to -q$ maps
strong stable manifolds homeomorphically
onto strong unstable ones and
therefore we obtain a family $\nu^{ss}$ of locally
finite Borel measures on strong stable manifolds
by defining $\nu^{ss}=\nu^{su}\circ {\cal F}$.
Let $dt$ be the usual Lebesgue measure on the
flow lines of the Teichm\"uller flow.
The locally finite Borel measure
$\tilde \nu$ on ${\cal Q}^1(S)$ defined
by $d\tilde \nu=d\nu^{ss}\times d\nu^{su}\times dt$
is invariant
under the Teichm\"uller geodesic flow $\Phi^t$ and the
action of the mapping class group.
If we denote by $\Delta$ the diagonal in
${\cal P\cal M\cal L}\times {\cal P\cal M\cal L}$ then the
desintegration $\hat \nu$ of $\tilde \nu$
along the flow lines of the Teichm\"uller flow
is an ${\cal M}(S)$-invariant locally finite
Borel measure on ${\cal P\cal M\cal L}\times
{\cal P\cal M\cal L}-\Delta$.
Let $\nu$ be the $\Phi^t$-invariant locally finite
Borel measure on ${\cal Q}(S)$ which
is the projection of the restriction of $\tilde \nu$
to a Borel fundamental domain for the action of
${\cal M}(S)$.
For the conformal density $\{\lambda^x\}$ in the
Lebesgue measure class the resulting
$\Phi^t$-invariant measure $\lambda$ on
${\cal Q}(S)$ is finite.

Call a quadratic differential $q\in {\cal Q}(S)$
\emph{forward returning} if there is a compact
subset $K$ of ${\cal Q}(S)$
and for every $k>0$ there is some $t>k$ with
$\Phi^tq\in K$.
Call a projective
measured geodesic lamination $\xi\in {\cal P\cal M\cal L}$
\emph{returning} if there is a quadratic
differential $q\in {\cal Q}^1(S)$
whose horizontal measured geodesic
lamination is contained in the class
$\xi$ and whose
projection to ${\cal Q}(S)$ is forward returning.
The set of returning projective measured geodesic
laminations is a Borel subset of
${\cal P\cal M\cal L}$
which is invariant
under the action of the mapping
class group and is contained
in the set of all uniquely
ergodic projective measured geodesic laminations
which fill up $S$
\cite{M82a}. Note however that Cheung and Masur
\cite{CM06} constructed an example of a
uniquely ergodic projective measured geodesic
lamination which fills up $S$ and is not returning.

Call a quadratic differential $q\in {\cal Q}(S)$
\emph{forward recurrent} if $q$ is contained in the
$\omega$-limit set of its own orbit under $\Phi^t$.
Call a projective measured geodesic lamination
$\xi\in {\cal P\cal M\cal L}$ \emph{recurrent} if
there is a quadratic differential $q\in {\cal Q}^1(S)$
whose horizontal measured geodesic lamination
is contained in the class $\xi$ and whose projection
to ${\cal Q}(S)$ is forward returning and contains 
a forward recurrent point $q_0\in {\cal Q}(S)$ in its
$\omega$-limit set.
The set ${\cal R\cal M\cal L}$ of all recurrent projective
measured geodesic laminations is an ${\cal M}(S)$-invariant
Borel subset of ${\cal P\cal M\cal L}$ which has full
Lebesgue measure.
We have.

\begin{lemma}\label{ergodicity}
For an ergodic
conformal density $\{\nu^x\}$ of dimension
$\alpha$ the following are equivalent.
\begin{enumerate}
\item $\{\nu^x\}$
gives full mass to the returning
projective measured geodesic laminations.
\item $\{\nu^x\}$ gives full mass to the set
${\cal R\cal M\cal L}$ of
recurrent projective measured geodesic laminations.
\item The measure $\hat\nu$ on
${\cal P\cal M\cal L}\times
{\cal P\cal M\cal L}-\Delta$ is ergodic under the
diagonal action of ${\cal M}(S)$.
\item The Teichm\"uller geodesic flow $\Phi^t$
is conservative for the measure $\nu$ on ${\cal Q}(S)$.
\item The Teichm\"uller geodesic flow $\Phi^t$
is ergodic for $\nu$.
\end{enumerate}
\end{lemma}
\begin{proof}
We follow Sullivan \cite{Su79} closely. Namely,
for $\epsilon >0$ let ${\rm int}{\cal T}(S)_{\epsilon}\subset {\cal T}(S)$
be the open ${\cal M}(S)$-invariant
subset of all complete hyperbolic structures
on $S$ of finite volume whose systole (i.e. the length
of a shortest closed geodesic) is bigger than
$\epsilon$ and define
${\cal Q}^1(\epsilon)=\{q\in
{\cal Q}^1(S)\mid Pq\in {\rm int}{\cal T}(S)_{\epsilon}\}$.
Then the closure $\overline{\cal Q}^1(\epsilon)$
of ${\cal Q}^1(\epsilon)$ projects to a compact
subset $\overline{\cal Q}(\epsilon)$ of ${\cal Q}(S)$.
For $\delta <\epsilon$ we have $\overline{\cal Q}(\delta)\supset
\overline{\cal Q}(\epsilon)$ and
$\cup_{\epsilon >0}\overline{\cal Q}(\epsilon)={\cal Q}(S)$.
For $\epsilon>0$ define ${\cal B}_\epsilon\subset
{\cal P\cal M\cal L}$ to be the
set of all projective measured geodesic laminations $\xi$
such that there is a quadratic differential
$q\in {\cal Q}^1(S)$ with $\pi(q)=\xi$ and a sequence
$\{t_i\}\to \infty$ with $\Phi^{t_i}q\in
{\cal Q}^1(\epsilon)$ for all $i>0$. The
set ${\cal B}_\epsilon$ is invariant
under the action of ${\cal M}(S)$, and
${\cal B}=\cup_{\epsilon >0}{\cal B}_\epsilon$
is the set of all returning projective measured geodesic
laminations.

Let $\{\nu^x\}$ be an ergodic conformal
density of dimension $\alpha$ which gives full
mass to the set ${\cal B}$ of returning projective
measured geodesic laminations.
By invariance and
ergodicity, there is a number $\epsilon >0$ such that
$\nu^x$ gives full mass to ${\cal B}_\epsilon$.
Let $\nu$ be the $\Phi^t$-invariant Radon measure
on ${\cal Q}(S)$ defined by $\{\nu^x\}$.
By the results of \cite{M80,M82a}, for
every forward returning quadratic differential
$q\in {\cal Q}(S)$ and every $z\in W^{ss}(q)$ we have
$d(\Phi^tq,\Phi^tz)\to 0$ $(t\to \infty)$. Thus
for $\nu$-almost every $q\in {\cal Q}(S)$
the $\Phi^t$-orbit of $q$
enters the compact set
$\overline{\cal Q}(\epsilon/2)$ for arbitrarily large times.
Then there is a number $\delta >0$ and for
$\nu$-almost every $q\in \overline{\cal Q}(\delta)$
there are infinitely many integers $m>0$ with
$\Phi^mq\in \overline{\cal Q}(\delta)$.
As a consequence, the first return map
to $\overline{\cal Q}(\delta)$
of the homeomorphism $\Phi^1$ of ${\cal Q}(S)$
defines
a measurable map $G:\overline{\cal Q}(\delta)\to
\overline{\cal Q}(\delta)$ which preserves
the restriction $\nu_0=\nu{\vert \overline{\cal Q}(\delta)}$
of $\nu$. Since $\nu$
is a Radon measure, $\nu_0$ is finite and hence
the system $(\overline{\cal Q}(\delta),\nu_0,G)$
is conservative. But then the measure $\nu$ is
conservative for the time-one map $\Phi^1$ of the
Teichm\"uller flow. Moreover, by the Poincar\'e recurrence
theorem, applied to the measure preserving map
$G:\overline{\cal Q}(\delta)\to \overline{\cal Q}(\delta)$,
we obtain that $\nu$-almost every $q\in {\cal Q}(S)$ is forward
recurrent. Thus 1) above implies 2) and 4).

Using the usual Hopf argument \cite{Su79} we conclude
that the measure $\nu$ is ergodic. Namely,
choose a continuous positive function $\rho:{\cal Q}(S)\to
(0,\infty)$ such that $\int \rho d\nu=1$; such a
function $\rho$ exists since the measure $\nu$ is locally
finite by assumption. By the above,
if $q,q^\prime\in {\cal Q}(S)$ are typical for $\nu$
and contained in the same strong
stable manifold then the orbits of $\Phi^t$ through
$q,q^\prime$ are forward asymptotic \cite{M80,M82a}.
By the Birkhoff ergodic theorem, for every
continuous function $f$ with compact support
the limit $\lim_{t\to \infty}\int_0^Tf(\Phi^t q)dt/
\int_0^T\rho(\Phi^t q)dt=f_\rho(q)$
exists almost everywhere, and the Hopf argument shows
that the function $f_\rho$
is constant along $\nu$-almost
all strong stable and strong unstable manifold.
From this ergodicity follows as in \cite{Su79}.
In particular, 1) above implies 5).

The remaining implications in the statement
of the lemma are either trivial or standard.
Since they will not be used in the sequel,
we omit the proof.
\end{proof}

Every ergodic conformal density either gives full measure
or zero measure
to the ${\cal M\cal L}$-invariant
Borel subset ${\cal R\cal M\cal L}$ of recurrent
projective measured geodesic laminations. The main goal of this
section is to show that a conformal density $\{\nu^x\}$
which gives
full measure to ${\cal R\cal M\cal L}$ is contained
Lebesgue measure class. For this we adapt ideas of Sullivan
\cite{Su79} to our situation. Namely, the set
${\cal R\cal M\cal L}$ can be viewed as the radial limit set
for the action of ${\cal M}(S)$ on ${\cal T}(S)$, where
${\cal P\cal M\cal L}$ is identified with the Thurston boundary
of ${\cal T}(S)$. This means that every point in
${\cal R\cal M\cal L}$ is contained in a nested sequence of
neighborhoods which are images of a fixed set under elements
of the mapping class group. The mass of these sets with
respect to the measures $\nu^x,\lambda^x$ (where as before,
$\{\lambda^x\}$ is a conformal density in the Lebesgue measure
class) are controlled, and this
allows a comparison of measures.

To carry out this approach, we have to construct such
nested sequences of neighborhoods of points in ${\cal R\cal M\cal L}$
explicitly. We use the Gromov distances on the boundary
$\partial {\cal C}(S)$ of the curve graph ${\cal C}(S)$
for this purpose. This
boundary consists of all
unmeasured minimal geodesic laminations which fill up
$S$, equipped with a coarse Hausdorff
topology. If we denote by
${\cal F\cal M\cal L}\subset {\cal P\cal M\cal L}$
the Borel set of all projective measured geodesic
laminations whose support is a geodesic lamination
which is minimal and fills up $S$ then the natural
forgetful map
\begin{equation}
F:{\cal F\cal M\cal L}\to \partial {\cal C}(S)
\end{equation}
which assigns to a projective measured geodesic
lamination in ${\cal F\cal M\cal L}$ its support
is a continuous ${\cal M}(S)$-equivariant
surjection \cite{Kl99,H06a}. The
restriction of the projection map $F$
to the Borel set ${\cal R\cal M\cal L}
\subset {\cal F\cal M\cal L}$ is injective \cite{M82a}.

For $q\in {\cal Q}^1(S)$ and $r>0$ define
\begin{equation}\label{balldefinition}
B(q,r)=\{u\in W^{su}(q)\mid d(Pq,Pu)\leq r\}.
\end{equation}
Then $B(q,r)$ is a compact neighborhood of
$q$ in $W^{su}(q)$
with dense interior which depends
continuously on $q$ in the following sense.
If $q_i\to q$ in ${\cal Q}^1(S)$ then
$B(q_i,r)\to B(q,r)$ in the Hausdorff topology
for compact subsets of ${\cal Q}^1(S)$.
We have.

\begin{lemma}\label{closed}
\begin{enumerate}
\item
The map $F:{\cal F\cal M\cal L}\to \partial{\cal C}(S)$ is
continuous and closed.
\item If the horizontal measured geodesic
lamination of the quadratic differential $q\in {\cal Q}^1(S)$
is uniquely ergodic then the sets
$F(\pi(B(q,r))\cap {\cal F\cal M\cal L})$
$(r>0)$ form a neighborhood basis for
$F(\pi(q))$ in $\partial{\cal C}(S)$.
\end{enumerate}
\end{lemma}
\begin{proof} The first part of the lemma is
immediate from the description of the Gromov
boundary of ${\cal C}(S)$ in \cite{Kl99,H06a}.

To show the second part, note first that
for every $q\in {\cal Q}^1(S)$ the restriction of 
the projection
$\pi:W^{su}(q)\to {\cal P\cal M\cal L}$ is a homeomorphism
of $W^{su}(q)$ onto an \emph{open} subset of 
${\cal P\cal M\cal L}$. To see this, identify
${\cal P\cal M\cal L}$ with a section
$\Sigma$ of the fibration ${\cal M\cal L}\to {\cal P\cal M\cal L}$.
Then $\xi\in \Sigma$ is contained in $\pi W^{su}(q)$ if and only if
the function
$\zeta\to i(\xi,\zeta)+i(\pi(-q),\zeta)$ on $\Sigma$ 
is \emph{positive}. Since the intersection form $i$ on 
${\cal M\cal L}$ is continuous and since
$\Sigma$ is compact, this is an open condition for
points $\xi\in \Sigma$.

Now if $q\in {\cal Q}^1(S)$ is uniquely ergodic then
by \cite{Kl99} we have
$F^{-1}(F(\pi(q)))=\{\pi(q)\}$ and hence if
$r>0$ is arbitrary then
$F({\cal F\cal M\cal L}-{\rm int}B(q,r))$
is a closed subset of $\partial{\cal C}(S)$ which
does not contain $F(\pi(q))$. In particular,
\[F(\pi (B(q,r))\cap {\cal F\cal M\cal L})\]
is a neighborhood of $F(\pi(q))$ in
$\partial{\cal C}(S)$.
From this and continuity of
$F$ the second part of the lemma follows.
\end{proof}

Define
\begin{equation}\label{cala}
{\cal A}=\pi^{-1}({\cal F\cal M\cal L})\subset {\cal Q}^1(S).
\end{equation}
Recall from Section 2 that
there is a number $\kappa >0$ and there is
an ${\cal M}(S)$-equivariant family of distance
functions $\delta_h$ ($h\in
{\cal T}(S)$)
on $\partial {\cal C}(S)$
such that $\delta_h\leq e^{\kappa d(h,z)}
\delta_z$ for all $h,z\in {\cal T}(S)$.
For $q\in {\cal A}$ and
$\chi >0$ define $D(q,\chi)\subset \partial
{\cal C}(S)$ to be the
closed $\delta_{Pq}$-ball of radius $\chi$
about $F\pi(q)\in \partial {\cal C}(S)$.
The following lemma is a translation
of hyperbolicity of the curve graph into
properties of the distance functions $\delta_h$.

\begin{lemma}\label{ballestimate}
\begin{enumerate}
\item For every $\beta>0$
there is a number
$\rho=\rho(\beta) >0$ such that
\[D(\Phi^tq,\rho)\subset D(q,\beta)\]
for every $q\in {\cal A}$ and every
$t\geq 0$.
\item There is a number $\beta_0>0$ with the
following property. For every $q\in {\cal A}$ and every
$\epsilon >0$ there is a number $T(q,\epsilon)>0$
such that $D(\Phi^tq,\beta_0)\subset D(q,\epsilon)$ for
every $t\geq T(q,\epsilon)$.
\end{enumerate}
\end{lemma}
\begin{proof} By the results of \cite{MM99} 
(see Theorem 4.1 of
\cite{H07a} for an explicit statement), there is a number $L>0$
such that the image under $\Upsilon_{\cal T}$
of every Teichm\"uller geodesic is an
\emph{unparametrized $L$-quasi-geodesic in
${\cal C}(S)$}.
This means that for every
$q\in {\cal Q}^1(S)$ there is an increasing
homeomorphism $\sigma_q:\mathbb{R}\to \sigma_q(\mathbb{R})
\subset \mathbb{R}$ such that
the curve $t\to \Upsilon_{\cal T}(P\Phi^{\sigma_q(t)}q)$ is
an $L$-quasi-geodesic in ${\cal C}(S)$.

If $q\in {\cal A}$
then we have $\sigma_q(t)\to \infty$ $(t\to \infty)$ and
the unparametrized $L$-quasi-geodesic 
$t\to \Upsilon_{\cal T}(P\Phi^tq)$ converges as
$t\to \infty$ in ${\cal C}(S)\cup \partial{\cal C}(S)$
to the point $F(\pi(q))\in \partial{\cal C}(S)$
(see \cite{Kl99,H05,H06a}).
In particular, for $q\in {\cal A}$ and every $T>0$ there
is a number $\tau=\tau(q,T)>0$
such that
$d(\Upsilon_{\cal T}(P\Phi^{t}q),\Upsilon_{\cal T}(Pq))\geq T$
for all $t\geq \tau$.

Since ${\cal C}(S)$ is 
a hyperbolic geodesic metric space, any finite subarc of
an $L$-quasi-geodesic is contained in a tubular neighborhood
of a geodesic in ${\cal C}(S)$ 
of uniformly bounded radius. This implies that 
there is no backtracking along an unparametrized $L$-quasi-geodesic:
There is a constant
$b>0$ only depending on $L$ and the hyperbolicity
constant of ${\cal C}(S)$ such that
$d(\gamma(t),\gamma(0))\geq d(\gamma(s),\gamma(0))$ for
all $t\geq s\geq 0$ and every $L$-quasi-geodesic
$\gamma:[0,\infty)\to {\cal C}(S)$.
From the definition of the Gromov distances
$\delta_c$ $(c\in {\cal C}(S))$ we obtain that
there is a number $\alpha>0$
such that for every $L$-quasi-geodesic
ray $\gamma:[0,\infty)\to {\cal C}(S)$ with endpoint $\xi\in
\partial{\cal C}(S)$ and every $t>0$ the Gromov distances
$\delta_{\gamma(t)}$ on $\partial {\cal C}(S)$ satisfy
\begin{equation}\label{deltagamma}
\delta_{\gamma(t)}\geq \alpha e^{\alpha d(\gamma(t),\gamma(0))}
\delta_{\gamma(0)}\end{equation}
on the $\delta_{\gamma(t)}$-ball of radius $\alpha$ about $\xi$.
Let $\kappa >0$ be as in inequality (\ref{delta}) from Section 2
and define $\beta_0=\alpha/\kappa^2$.

By inequality (\ref{deltacomparison}) from Section 2, for $q\in {\cal A}$
and $t\geq 0$ we have 
\begin{equation}\label{betanull}
\delta_{P\Phi^tq}\geq \kappa^{-2}\alpha 
e^{\alpha d(\Upsilon_{\cal T}(P\Phi^tq),\Upsilon_{\cal T}(Pq))}\delta_{Pq}
\end{equation}
on the $\delta_{P\Phi^tq}$-ball $D(\Phi^tq,\beta_0)$.
Thus if for $\epsilon >0$ we choose 
$T_1>0$ sufficiently large that 
$\epsilon\alpha e^{\alpha T_1}\geq\kappa^2\beta_0$ then
for $q\in {\cal A}$, for $T=\tau(q,T_1)>0$ and 
for $t>T$ we have $D(\Phi^t q,\beta_0)\subset D(q,\epsilon)$
which shows the second part of the lemma. 

To show the first part of the lemma, for $\beta <\beta_0$ define
$\rho(\beta)=\alpha\beta/\kappa^2$. Then the
estimate (\ref{betanull}) above shows that
$D(\Phi^tq,\rho)\subset D(q,\beta)$ for every $q\in {\cal A}$ and
every $t\geq 0$. This completes the proof of the lemma.
\end{proof}

For a forward recurrent point
$q_0\in {\cal Q}(S)$ let
${\cal R\cal M\cal L}(q_0)\subset {\cal R\cal M\cal L}$
be the Borel subset of all recurrent projective
measured geodesic laminations $\xi\in {\cal R\cal M\cal L}$
such that there is some $q\in \pi^{-1}(\xi)$
with the following property. The
projection to ${\cal Q}(S)$ of the orbit of $q$ under
the Teichm\"uller geodesic flow contains
$q_0$ in its $\omega$-limit set.
By definition, the set ${\cal R\cal M\cal L}(q_0)$ is
invariant under the action of ${\cal M}(S)$. Moreover,
every recurrent point $\xi\in {\cal R\cal M\cal L}$
is contained in one of the sets ${\cal R\cal M\cal L}(q_0)$
for some forward recurrent point $q_0\in {\cal Q}(S)$.
Note moreover that an orbit of $\Phi^t$
in ${\cal Q}(S)$ which is typical for the $\Phi^t$-invariant
Lebesgue measure on ${\cal Q}(S)$ is dense and hence
for every forward recurrent point
$q_0\in {\cal Q}(S)$ the
set ${\cal R\cal M\cal L}(q_0)$ has full Lebesgue
measure.
Write
\begin{equation}
C(q_0)=F({\cal R\cal M\cal L}(q_0))\subset
\partial{\cal C}(S)\quad\text{and}\quad
A(q_0)=\pi^{-1}{\cal R\cal M\cal L}(q_0)\subset {\cal Q}^1(S).
\end{equation}

Following \cite{F69}, a \emph{Borel covering
relation} for a Borel subset $C$ of a metric
space $(X,d)$ is a family ${\cal V}$ of
pairs $(x,V)$ where $V\subset X$ is a Borel set,
where $x\in V$ and such that
\begin{equation}
C\subset \bigcup \{V\mid (z,V)\in {\cal V}\text{ for some }
z\in C\}.\end{equation} The covering
relation ${\cal V}$ is called
\emph{fine} at every point of $C$ if
for every $x\in C$ and every $\alpha >0$
there is some $(y,V)\in {\cal V}$ with
$x\in V\subset U(x,\alpha)$ where $U(x,\alpha)$ denotes
the open ball of radius $\alpha$ about $x$.

For $\chi>0$ and the forward recurrent point $q_0\in {\cal Q}(S)$
with lift $q_1\in {\cal Q}^1(S)$ there
is by continuity of $F\circ\pi$ a compact neighborhood $K$
of $q_1$ in ${\cal Q}^1(S)$ such that $F\circ\pi(K\cap {\cal A})\subset
D(q_1,\chi)$. We call $K$ a \emph{$\chi$-admissible} neighborhood
of $q_1$. For a number $\chi>0$ and such a $\chi$-admissible
neighborhood $K$ of $q_1$ define
\begin{align}\label{vitalirelation}
{\cal V}_{q_0,\chi,K}=& \{(F\pi(q),g D(q_1,\chi))\mid \\
& q\in W^{su}(q_1)\cap A(q_0),g\in {\cal M}(S),gK\cap \cup_{t>0}\Phi^tq\not=
\emptyset.\}\notag
\end{align}
We sometimes identify a pair
$(\xi,gD(q_1,\chi))\in {\cal V}_{q_0,\chi,K}$ with the set $gD(q_1,\chi)$
whenever the point $\xi$ has no importance.
Let $\beta_0>0$ be as in Lemma \ref{ballestimate}.
We have.

\begin{lemma}\label{fine}
Let $q_0\in {\cal Q}(S)$ be a forward recurrent
point and let $q_1\in {\cal Q}^1(S)$
be a lift of $q_0$. Then for every
$\chi<\beta_0/4$ and every $\chi$-admissible
compact neighborhood $K$ of $q_1$ the family
${\cal V}_{q_0,\chi,K}$
is a Borel covering relation
for $C(q_0)\subset (\partial {\cal C}(S),\delta_{Pq_1})$ which is fine
at every point of $C(q_0)$.
\end{lemma}
\begin{proof}
Using the above notations,
let $q_0\in {\cal Q}(S)$ be a forward
recurrent point
and let $q_1$ be a lift
of $q_0$ to ${\cal Q}^1(S)$.
Let $\chi<\beta_0/4$ where $\beta_0>0$ is as in
Lemma \ref{ballestimate}. It clearly suffices to show
the lemma for the covering relations ${\cal V}_{q_0,\chi,K}$
where $K$ is a sufficiently small
$\chi$-admissible neighborhood of $q_1$.

By relation (\ref{delta})
in Section 2 we infer
that for every sufficiently small $\chi$-admissible
neighborhood $K$ of $q_1$ we have
\begin{equation}
\delta_{Pq_1}/2\leq \delta_{Pq}\leq 2\delta_{Pq_1}\text { for
every }q\in K. \end{equation}
In particular, for $q\in K\cap {\cal A}$ the
set $D(q_1,\chi)$ contains $F\pi(q)$ and
is contained in $D(q,4\chi)$.

By the construction of the distances
$\delta_h$ $(h\in {\cal T}(S))$ on $\partial {\cal C}(S)$
it suffices to show that for every $q\in A(q_0)$ and every
$\epsilon>0$ there is some $g\in {\cal M}(S)$ such that
the set $gD(q_1,\chi)$ contains $F(\pi(q))$ and is contained
in the open $\delta_{Pq}$-ball
of radius $\epsilon$ about $F(\pi(q))$.

For $q\in A(q_0)$
and $\epsilon >0$ let $T(q,\epsilon)>0$ be as in the second
part of Lemma \ref{ballestimate}.
Choose some $t>T(q,\epsilon)$ such that
$\Phi^tq\in \tilde K=\cup_{g\in {\cal M}(S)}K$;
such a number exists by
the definition of the set $A(q_0)$ and
by \cite{M80}. By Lemma \ref{ballestimate}
we have $D(\Phi^tq,4\chi)\subset D(q,\epsilon)$.
Now if
$g\in {\cal M}(S)$ is such that $\Phi^tq\in gK$
then we obtain from $\chi$-admissibility of
the set $K$
and equivariance
under the action of the mapping class group that
\begin{equation}
F\pi(q)=F\pi(\Phi^tq)\in gD(q_1,\chi)\subset
D(\Phi^tq,4\chi)\subset D(q,\epsilon).
\end{equation} Since
$\epsilon >0$ was arbitrary, this shows the lemma.
\end{proof}

The next proposition is the main technical result
of this section. For its formulation, recall from
\cite{F69} the definition of a \emph{Vitali relation}
for a finite Borel measure on the Borel subset $C(q_0)$
of ${\cal C}(S)$. We show.

\begin{proposition}\label{Vitali}
Let $q_0\in {\cal Q}(S)$ be a forward recurrent point
for the Teichm\"uller geodesic flow. Then for
every sufficiently small $\chi>0$, every
sufficiently small $\chi$-admissible compact neighborhood
$K$ of $q_1$ and for every
conformal density $\{\nu^x\}$
on ${\cal P\cal M\cal L}$ which
gives full measure to the set ${\cal R\cal M\cal L}(q_0)$,
the covering relation
${\cal V}_{q_0,\chi,K}$ for $C(q_0)$ is a Vitali relation for
the measure $F_*\nu^x$ on $\partial{\cal C}(S)$.
\end{proposition}
\begin{proof}
The strategy of proof is to use the properties of the
balls $D(q,\epsilon)$ established in Lemma \ref{ballestimate}
to gain enough control on
$\nu^x$-volumes that the results of Federer \cite{F69}
can be applied.

Let $q_0\in {\cal Q}(S)$
be a forward recurrent point for
$\Phi^t$ and let
$q_1\in {\cal Q}^1(S)$ be a lift of $q_0$.
Since no torsion element of ${\cal M}(S)$
fixes pointwise the Teichm\"uller geodesic defined
by $q_1$
we may assume that
the point $Pq_1\in {\cal T}(S)$ is
not fixed by any nontrivial element of the mapping
class group.

By Lemma \ref{fine} and using the notations
from this lemma, for every
$\chi<\beta_0/4$ and every $\chi$-admissible
compact neighborhood $K$ of $q_1$ the covering
relation ${\cal V}_{q_0,\chi,K}$ for
$C=F({\cal R\cal M\cal L}(q_0)-\pi(-q_1))
\subset \partial {\cal C}(S)$ is fine
for the metric $\delta_{Pq_1}$ at every point of
$C$.

We first establish some geometric control
on the covering relation ${\cal V}_{q_0,\chi,K}$
for some particularly chosen small $\chi<\beta_0/4$ and
a suitably chosen $\chi$-admissible neighborhood $K$
of $q_1$. For this
let again $d$ be the distance on ${\cal T}(S)$
defined by the Teichm\"uller metric.
Choose a number $r>0$ which is sufficiently small
that the images
under the action of the mapping class
group of the closed $d$-ball $B(Pq_1,5r)$ of radius $5r$
about $Pq_1$ are pairwise disjoint. By the estimate
(\ref{delta}) for the family of distance functions
$\delta_z$ $(z\in {\cal T}(S))$, by decreasing the size of
the radius $r$ we may assume that
\begin{equation}\label{deltaxu}
\delta_x/2\leq \delta_u\leq 2\delta_x\text{ for all }x,u\in B(Pq_1,5r).
\end{equation}

Recall from (\ref{balldefinition})
above the definition of the closed
balls $B(q,r)\subset W^{su}(q)$ $(q\in {\cal Q}^1(S))$.
By continuity of the projection $\pi$ there
is an open neighborhood $U_1$ of $q_1$ in
${\cal Q}^1(S)$ such that
\begin{equation}
\pi B(z,r)\subset \pi B(q,2r)
\text{ for all }z,q\in U_1
\end{equation} and that moreover
the projection $PU_1$ of $U_1$ to ${\cal T}(S)$ is
contained in the open ball of radius $r$ about
$Pq_1$. This implies in particular that
$gU_1\cap U_1=\emptyset$ for every nontrivial
element $g\in {\cal M}(S)$.

Since the horizontal measured geodesic lamination
of $q_1$ is uniquely ergodic, Lemma \ref{closed} shows
that
there is a number $\beta>0$ such that
\begin{equation}
F(\pi B(q_1,r)\cap {\cal F\cal M\cal L})\supset
D(q_1,8\beta).
\end{equation}

Since the map $F\circ \pi:{\cal A}
\to \partial {\cal C}(S)$ is continuous, there is
an open neighborhood $U_2\subset U_1$ of $q_1$
such that
$U_2\cap {\cal A}\subset (F\circ\pi)^{-1}D(q_1,\beta)$.
By the choice of $U_1$ we have
$D(z,\beta)\subset D(q,8\beta)$ for all
$q,z\in U_2\cap {\cal A}$.
For all $q,z\in U_2\cap {\cal A}$
we also have
\begin{equation}\label{eins}
D(z,\beta)\subset D(q_1,4\beta)\subset
F(\pi B(q_1,r)\cap {\cal F\cal M\cal L})\subset
F(\pi B(q,2r)\cap {\cal F\cal M\cal L}).
\end{equation}

By Lemma \ref{ballestimate},
there is a number $\sigma\leq \beta$
such that for every
$t\geq 0$ we have
\begin{equation}
D(\Phi^tq,\sigma)\subset D(q,\beta).
\end{equation}

Now $U_2$ is an open neighborhood of $q_1$ in
${\cal Q}^1(S)$ and therefore
$U_2\cap W^{su}(q_1)$ is an open neighborhood
of $q_1$ in $W^{su}(q_1)$. In particular,
there is a number $r_1<r$ such that
$B(q_1,r_1)\subset W^{su}(q_1)\cap U_2$. Thus
by Lemma \ref{closed}
there is a number $\chi\leq \sigma/16$ such that
\begin{equation}\label{zwei}
F(\pi(W^{su}(q_1)\cap U_2)\cap {\cal F\cal M\cal L})
\supset D(q_1,16\chi).\end{equation}
Note that we have
\begin{equation}\label{zweib}D(\Phi^tq,16\chi)
\subset D(q,\beta)
\text{ for all}\,\, q\in U_2\cap {\cal A}
\,\,\text{and
all}\,\, t>0.
\end{equation}

Using once more continuity of the map
$F\circ\pi:{\cal A}\to
\partial {\cal C}(S)$,
there is a compact neighborhood $K\subset U_2$
of $q_1$ such that
\begin{equation}
K\cap {\cal A}\subset(F\circ\pi)^{-1}D(q_1,\chi).
\end{equation}
In particular, $K$ is $\chi$-admissible.
By inequality (\ref{deltaxu}) for
the dependence
of the metrics $\delta_{Pq}$ on the
points $q\in K\subset U_1$ we then have
\begin{equation}\label{sixteen}
F\pi(z)\in D(q_1,\chi)\subset
D(z,4\chi)\subset D(q_1,16\chi)
\text{ for all }z\in K\cap {\cal A}.
\end{equation}
By (\ref{zwei})
above and continuity of the strong
unstable foliation and of the map $\pi$
we may moreover assume that
\begin{equation}\label{drei}
F(\pi (W^{su}(q)\cap U_2)\cap {\cal F\cal M\cal L})
\supset D(q_1,16\chi)
\text{ for every }q\in K.
\end{equation}
Note that
if $z\in K\cap {\cal A}$ then $D(z,4\chi)\subset D(q_1,16\chi)$
and hence if $u\in W^{su}(q_1)\cap A(q_0)$
is such that $F\pi(u)\in D(z,4\chi)$ then
$u\in U_2$. Namely, the projective measured geodesic
lamination $\pi(u)$ of
every $u\in A(q_0)$ is uniquely ergodic
and therefore $(F\circ \pi)^{-1}(F(\pi(u))\cap W^{su}(q_1)$
consists of a unique point. However, by (\ref{zwei}) above,
the set $U_2\cap W^{su}(q_1)$ contains such a point.
Consequently,
inequality (\ref{deltaxu})
above shows that
$D(z,4\chi)\subset D(u,16\chi)$.

Define \begin{equation}\label{Vitali2}
{\cal V}_0={\cal V}_{q_0,\chi,K}.
\end{equation}
By Lemma \ref{fine},
${\cal V}_0$ is a covering relation for
the set $C\subset C(q_0)\subset \partial{\cal C}(S)$ which
is fine at every point of $C$.

By the choice of the set $K\subset U_1$,
if $q\in W^{su}(q_1)\cap A(q_0)$, if $g\in
{\cal M}(S)$ and if $t>0$ are such that
$\Phi^tq\in gK$ then $g\in {\cal M}(S)$ is
uniquely determined by $\Phi^tq$.
For $(\xi,V)\in {\cal V}_0$
define \begin{align}
\rho(\xi,V)=& \max\{e^{-t}\mid q\in W^{su}(q_1)\cap A(q_0),
t\geq 0,\\
& V=gD(q_1,\chi), \Phi^tq\in gK,\pi(q)=\xi\}.\notag
\end{align}

Following \cite{F69},
for $(\xi,V)\in {\cal V}_0$
define the \emph{$\rho$-enlargement}
of $V$ by \begin{equation}\label{vier}
\hat V=\bigcup \{W\mid (\zeta,W)\in {\cal V}_0,
W\cap V\cap C(q_0)\not=\emptyset, \rho(\zeta,W)\leq e^{r}\rho(\xi,V)\}
\end{equation}
where in this definition, the constant $r>0$ is chosen
as in the beginning of this proof.

Let $\{\nu^x\}$ be a conformal
density of dimension $\alpha\geq 0$
which gives full measure to the set ${\cal R\cal M\cal L}(q_0)$.
We may assume that
the density is ergodic, i.e. that the
${\cal M}(S)$-invariant
measure class it defines on ${\cal P\cal M\cal L}$
is ergodic. The measure $\nu^x$
induces a Borel measure $F_*\nu^x$ on
the set $C=F({\cal R\cal M\cal L}(q_0)-\pi(-q_1))
\subset C(q_0)\subset \partial {\cal C}(S)$.

Recall from the beginning of this section
that the measures $\nu^y$
$(y\in {\cal T}(S))$ define
a family of ${\cal M}(S)$-invariant
Radon measures $\nu^{su}$
on strong unstable
manifolds in ${\cal Q}^1(S)$. These measures are invariant
under holonomy along
strong stable manifolds and they are
quasi-invariant under the
Teichm\"uller geodesic flow, with transformation
$d\nu^{su}\circ \Phi^t=e^{\alpha t}d\nu^{su}$.
For $q\in {\cal Q}^1(S)$ the measure
$\nu^{su}$ on $W^{su}(q)$
projects to a Borel measure $\nu_q$ on $C$.
For $q,z\in {\cal Q}^1(S)$
the measures $\nu_q,\nu_z$ are absolutely continuous,
with continuous Radon Nikodym derivative
depending continuously on $q,z$.
By invariance of the measures $\nu^{su}$ under
holonomy along strong stable manifolds and by the choice
of the point $q_1$ and the number $\chi>0$
there is a
number $a>0$ such that $1/a\geq \nu_qD(q_1,\chi)\geq a$
for all $q\in K$.

Write $\nu_1=\nu_{q_1}$;
we claim that there is a number $c>0$ such that
$\nu_1(\hat V)\leq c\nu_1(V)$
for all $(\xi,V)\in {\cal V}_0$.
For this let $(\xi,V)\in {\cal V}_0$ be arbitrary;
then there is some $q\in W^{su}(q_1)\cap A(q_0)$
with $F\pi(q)=\xi$ and there is a number
$t\geq 0$ and some $g\in {\cal M}(S)$
such that $\Phi^tq\in gK$ and that
$V=gD(q_1,\chi)$ and $\rho(\xi,V)=e^{-t}$.
By equivariance under the action of the mapping
class group and by the inclusion (\ref{eins}) above, we
have
\begin{equation}\label{fuenf}
V=gD(q_1,\chi)\subset F(\pi B(\Phi^tq,2r)\cap
{\cal F\cal M\cal L}).\end{equation}

Let $(\zeta,W)\in {\cal V}_0$ be such that
\begin{equation}
\rho(\xi,V)\leq\rho(\zeta,W)\leq e^{r}\rho(\xi,V)
\end{equation}
and that $W\cap V\cap C\not=\emptyset$.
Then there is a number
$\epsilon \in [0,r]$, a point $z\in W^{su}(q_1)$
such that $F\pi(z)=\zeta$ and some $h\in {\cal M}(S)$
such that $\Phi^{t-\epsilon}z\in hK$ and that
$W=hD(q_1,\chi)$,
$\rho(\zeta,W)=e^{\epsilon-t}$.
By equivariance under the action of
${\cal M}(S)$ and the inclusion (\ref{eins}) above,
we have
\begin{equation}W=hD(q_1,\chi)\subset
F(\pi B(\Phi^{t-\epsilon}z,2r)\cap {\cal F\cal M\cal L})
\end{equation}
and hence from the definition of the sets $B(q,R)$ and
the definition of the strong unstable manifolds we
conclude that
\begin{equation}\label{sechs}
W\subset
F(\pi B(\Phi^tz,4r)\cap {\cal F\cal M\cal L}).
\end{equation}
Since the restriction of the map $F$ to
${\cal R\cal M\cal L}$
is injective and the restriction of the map
$\pi$ to $W^{su}(\Phi^tq)$ is injective and
since $V\cap W
\cap C\not=\emptyset$ by assumption we obtain that
$B(\Phi^tz,4r)\cap B(\Phi^tq,2r)\not=\emptyset$.
As a consequence, the distance
in ${\cal T}(S)$ between the
points $P(\Phi^{t}z)$
and $P(\Phi^tq)$ is at most $6r$ and hence
the distance between $P(\Phi^{t-\epsilon}z)\in PhK=hPK$
and $P(\Phi^tq)\in PgK=gPK$ is at most $7r$.
On the other hand, since $K\subset U_1$, for
$u\not=v\in {\cal M}(S)$ the distance in ${\cal T}(S)$
between $uPK$ and $vPK$ is not smaller than $8r$.
Therefore we have $g=h$ and
$V=W$.
This shows that
\begin{equation}\label{sieben}
\nu_1(\bigcup\{W\mid (\zeta,W)\in {\cal V}_0,
\rho(\xi,V)\leq\rho(\zeta,W)\leq
e^{r}\rho(\xi,V),W\cap V\cap C\not=\emptyset\})
= \nu_1(V).\end{equation}

On the other hand, if $z\in W^{su}(q_1)\cap A(q_0)$, if
$s\geq 0$ and $h\in {\cal M}(S)$ are such
that $\Phi^sz\in hK$ and $hD(q_1,\chi)=W$
and if
$(F\pi(z),W)\in {\cal V}_0$ is such that
\begin{equation}
e^{-s}=\rho(F\pi(z),W)\leq \rho(\xi,V)
\end{equation}
and $V\cap W \cap C\not=\emptyset$
then $s\geq t$. By
the choice of the set $K$, equivariance under the
action of the mapping class group
and the inclusion (\ref{sixteen}) above,
we have
$W\subset D(\Phi^sz,4\chi)$ and
$V\subset D(\Phi^tq,4\chi)$ and hence
$D(\Phi^tq,4\chi)\cap
D(\Phi^sz,4\chi)\cap C\not=\emptyset$.
In other words,
there is some $u\in A(q_0)\cap
W^{su}(\Phi^tq)$ with
$F(\pi(u))\in D(\Phi^tq,4\chi)\cap D(\Phi^sz,4\chi)$.

By the inclusions (\ref{sixteen}) and (\ref{drei}) and the following
remark, since
$u\in W^{su}(\Phi^tq)\cap A(q_0),
\Phi^tq\in gK$ and $F\pi(u)\in
D(\Phi^tq,4\chi)\subset gD(q_1,\sigma)$
we have $u\in gU_2\cap A(q_0)$
and moreover
\begin{equation}\label{acht}
\Phi^{s-t}u\in
W^{su}(\Phi^sz)\cap A(q_0)\quad \text{and}\quad
W\subset D(\Phi^sz,4\chi)\subset D(\Phi^{s-t}u,16\chi).
\end{equation}
From (\ref{zweib}) above and invariance under the
action of the mapping class group we obtain
$D(\Phi^{s-t}u,16\chi)\subset D(u,\beta)$.
The inclusion (\ref{eins})
then yields that
\begin{equation}\label{neun}
W\subset D(\Phi^{s-t}u,16\chi)
\subset D(u,\beta)\subset
F(\pi B(\Phi^tq,2r)\cap {\cal F\cal M\cal L}).
\end{equation}
This shows that the $\rho$-enlargement
$\hat V$ of $V$ is contained in
$F(\pi B(\Phi^tq,2r)\cap {\cal F\cal M\cal L})$.

Since $\Phi^tq\in \cup_{g\in {\cal M}(S)}gK$ by assumption,
by invariance under the action of the mapping
class group
the $\nu_{\Phi^tq}$-mass of
$F(\pi B(\Phi^tq,2r)\cap {\cal F\cal M\cal L})$ is
uniformly bounded. Therefore by
the transformation rule for the measures
$\nu_z$ under the action of the Teichm\"uller
geodesic flow,
the $\nu_1$-mass of $\hat V$ is
bounded from above by a fixed multiple
of the $\nu_1$-mass of $V$.
Thus by the results of Federer \cite{F69} and
by Lemma \ref{fine},
the covering relation ${\cal V}_{q_0,\chi,K}$ is indeed a
Vitali relation for the measure
$F_*\nu_1$ and hence it is a Vitali relation
for the measure $F_*\nu^x$ as well. Note that the same is true
for the covering relation
${\cal V}_{q_0,\epsilon,K^\prime}$ for every $\epsilon <\chi$
and every sufficiently small $\epsilon$-admissible neighborhood $K^\prime$
of $q_1$.
\end{proof}

Using Lemma \ref{ergodicity}, Lemma \ref{fine} and
Proposition \ref{Vitali}
we can now show.

\begin{lemma}\label{confdim}
 \begin{enumerate}
\item A conformal density
on ${\cal P\cal M\cal L}$ which gives full
measure to ${\cal F\cal M\cal L}$ is of dimension
at least $6g-6+2m$, with equality if and only
if it coincides with the Lebesgue measure up to scale.
\item A conformal density
which gives full measure to the set of returning projective
measured geodesic laminations
is of dimension $6g-6+2m$ and coincides with
the Lebesgue measure up to scale.
\end{enumerate}
\end{lemma}
\begin{proof} Let $\{\nu^x\}$ be a conformal
density of dimension $\alpha\geq 0$ which gives
full measure to the set ${\cal F\cal M\cal L}$.
We may assume that
the density is ergodic, i.e. that the
${\cal M}(S)$-invariant
measure class it defines on ${\cal P\cal M\cal L}$
is ergodic. Let moreover $\{\lambda^x\}$ be
the conformal density of dimension
$h=6g-6+2m$ which defines the $\Phi^t$-invariant
probability measure on ${\cal Q}(S)$ in the Lebesgue
measure class.
We have to show that $\alpha\geq h$, with equality
if and only if $\nu^x=\lambda^x$ up to scale. For this assume
that $\alpha \leq h$.

The Lebesgue measure $\lambda$ on ${\cal Q}(S)$
is of full support and ergodic under the Teichm\"uller flow and therefore
the $\Phi^t$-orbit of $\lambda$-almost every $q\in {\cal Q}(S)$
is dense. This implies that
there is a recurrent point
$q_0\in {\cal Q}(S)$
with the property that the measure $\lambda^x$ gives
full mass to ${\cal R\cal M\cal L}(q_0)$.

Recall that the conformal
densities $\{\nu^x\},\{\lambda^x\}$ define
families $\nu^{su},\lambda^{su}$ of Radon measures
on the strong unstable manifolds. For $q\in {\cal Q}^1(S)$
denote by $\nu_q,\lambda_q$ the image under the
map $F\circ\pi$ of the
restriction of these measures
to $W^{su}(q)\cap {\cal A}$. Since the conformal densities
$\{\nu^x\},\{\lambda^x\}$ give full measure to the
set ${\cal F\cal M\cal L}$, for every $q\in {\cal Q}^1(S)$ the
measures $\lambda_q,\nu_q$ on $\partial{\cal C}(S)$ are of full
support.

Let $q_1\in {\cal Q}^1(S)$ be a lift of $q_0$ to
${\cal Q}^1(S)$.
By Proposition \ref{Vitali}, for sufficiently small
$\chi>0$ and for a sufficiently
small compact neighborhood $K$ of $q_1$
the covering relation
${\cal V}_{q_0,\chi,K}$ is a Vitali relation
for the measure
$\lambda_{q_1}$ on $\partial {\cal C}(S)$.
Using equivariance under the action of ${\cal M}(S)$
and the fact that
$\nu_{q_1}$ is of full support,
if the measures $\nu_{q_1},\lambda_{q_1}$ are singular
then there
for $\lambda^{su}$-almost every $q\in W^{su}(q_1)$ there
is a sequence
$t_i\to \infty$
such that for every $i>0$ the following holds.
\begin{enumerate}
\item $\Phi^{t_i}q\in g_iK$ for some $g_i\in {\cal M}(S)$.
\item The $\nu_{\Phi^{t_i}q}$-mass and the
$\lambda_{\Phi^{t_i}q}$-mass of $D(g_iq_1,\chi)$ is
bounded from above and below by a universal constant.
\item The limit
$\lim_{i\to \infty}\nu_{q_1}(D(g_iq_1,\chi))/
\lambda_{q_1}(D(g_iq_1,\chi))$ exists and equals zero.
\end{enumerate}
In particular, for every $k>0$ and
all sufficiently large
$i$, say for all $i\geq i(k)$,
we have $\lambda_{q_1}D(g_iq_1,\chi)\geq k
\nu_{q_1}D(g_iq_1,\chi)$.
On the other hand, we also have
\begin{equation}\lambda_{q_1}D(g_iq_1,\chi)
=e^{-h t_i}\lambda_{\Phi^{t_i}q}D(g_iq_1,\chi)
\leq ce^{-h t_i}\end{equation}
for a universal constant $c>0$ and
$\nu_{q_1}D(g_iq_1,\chi)\geq de^{-\alpha t_i}$ for
a universal constant $d>0$. If $k>0$ is sufficiently large
that $kd \geq 2c$ then for
$i\geq i(k)$ we obtain a contradiction.

In other words, if $\alpha\leq h$ then the measures
$\{\nu^x\}$ and $\{\lambda^x\}$ are absolutely continuous.
Moreover, they give full mass to the recurrent projective
measured geodesic laminations.
Then they define absolutely continuous
$\Phi^t$-invariant Radon measures $\nu,\lambda$ on
${\cal Q}(S)$ which are ergodic by Lemma \ref{ergodicity}.
As a consequence,
the measures $\nu,\lambda$ coincide up to scale and hence
the measures
$\{\nu^x\},\{\lambda^x\}$
coincide up to scale as well. This shows the first part of the lemma.

To show the second part of the lemma, assume
that $\alpha \geq h$ and that the conformal density
$\{\nu^x\}$ gives full measure to the subset of
${\cal P\cal M\cal L}$ of returning points.
By Lemma \ref{ergodicity}, $\{\nu^x\}$ gives full measure
to the set ${\cal R\cal M\cal L}$ of recurrent points.
Since $\{\nu^x\}$ is ergodic, there is a forward recurrent
quadratic differntial $q_0\in {\cal Q}(S)$ such that
$\{\nu^x\}$ gives full measure to the set ${\cal R\cal M\cal L}(q_0)$.
This implies that we can
exchange the roles of $\{\lambda^x\}$ and $\{\nu^x\}$
in the above argument and obtain
that $\{\nu^x\},\{\lambda^x\}$ coincide up to scale.
\end{proof}

The following proposition uses the results
of Minsky and Weiss
\cite{MW02} to show that up to scale,
the Lebesgue measure is the unique conformal
density on ${\cal P\cal M\cal L}$ which
gives full measure to the set ${\cal F\cal M\cal L}$
of filling projective measured geodesic laminations.

\begin{proposition}\label{confonfilling}
Let $\{\nu^x\}$ be a conformal density which
gives full measure to ${\cal F\cal M\cal L}$. Then
$\{\nu^x\}$ is the Lebesgue measure up to scale.
\end{proposition}
\begin{proof} We argue by contradiction
and we assume that there
is a conformal density $\{\nu^x\}$
which gives full measure to the set ${\cal F\cal M\cal L}$
of all projective measured
geodesic laminations whose support is minimal and fills up $S$
and which is singular to the Lebesgue measure.
Without loss of generality we can assume that $\{\nu^x\}$ is
ergodic.
By Lemma \ref{confdim},
the dimension $\alpha$ of $\{\nu^x\}$ is strictly
bigger than $h=6h-6+2m$ and the 
$\nu^x$-measure of the set of returning points vanishes.

Let $\nu$ be the locally finite Borel measure on ${\cal Q}(S)$ which
can be written in the form
$d\nu=d\nu^{su}\times d\lambda^s$ where $\lambda^s$ is the
family of Lebesgue measures on stable manifolds which transforms
under the Teichm\"uller flow $\Phi^t$ via
$d\lambda^s\circ\Phi^t=e^{-ht}d\lambda^s$.
The measure $\nu$ is quasi-invariant under the Teichm\"uller
geodesic flow and transforms via
\[\nu\circ \Phi^t=e^{(\alpha-h)t}\nu.\]
Since $\alpha >h$ this implies that the measure
$\nu$ is infinite.
Moreover, it gives full mass to quadratic differential
whose horizontal measured geodesic lamination
is minimal and fills up $S$.

The family $\lambda^s$ of Lebesgue measures on stable
manifolds is invariant under the horocycle flow
$h_t$ as defined in Section 2. By the explicit construction
of the measures $\nu^{su}$, this implies that the measure $\nu$ is
invariant under $h_t$.

However, following the reasoning
of Dani \cite{D79} (see the proof of
Corollary 2.6 of \cite{MW02} for a discussion in our
context), this implies that
the measure $\nu$ is necessarily \emph{finite}.
Namely, by the Birkhoff
ergodic theorem, applied to the horocycle flow $h_t$
and the locally finite $h_t$-invariant
measure $\nu$ (see Theorem 2.3 of
\cite{K85} for the version of the Birkhoff ergodic theorem
for locally finite infinite measures needed here),
for $\nu$-almost every $q\in {\cal Q}(S)$ and
for every continuous positive
function $f$ on ${\cal Q}(S)$ with $\int  f d\nu<\infty$
the limit \begin{equation}
\lim_{T\to \infty}\frac{1}{T}\int_0^T f(h_tq)dt=F(q)
\end{equation}
exists, and the resulting
function $F$ is $h_t$-invariant and $\nu$-integrable.
On the other hand, consider
the family of Borel probability measures
\begin{equation}
\mu(q,T)=\frac{1}{T}\int_0^T\delta(h_tq)dt
\end{equation}
on ${\cal Q}(S)$
where $\delta(x)$ denotes the Dirac mass at $x$.
By Theorem H2 of \cite{MW02} (more precisely, by the theorem
in the appendix which is
a slightly extended version of this result), for
every $\epsilon >0$ there is a compact set
$K_\epsilon\subset {\cal Q}(S)$
such that for $\nu$-almost every $q\in {\cal Q}(S)$,
every weak limit $\mu(q,\infty)$ of the measures
$\mu(q,T)$ as $T\to \infty$ satisfies
$\mu(q,\infty)(K_\epsilon)\geq 1-\epsilon$.
Since the function $f$ is \emph{positive},
we have $\inf\{f(z)\mid z\in K_{1/2}\}=2c >0$ and therefore
$F(q)=\int fd\mu(q,\infty)\geq  c$ for $\nu$-almost every
$q\in {\cal Q}(S)$. But this contradicts the fact that the
measure $\nu$ is infinite and that $F$ is $\nu$-integrable and
shows the proposition.
\end{proof}

\section{Train tracks}

In Section 3 we showed that conformal densities
for the mapping class group which give full
measure to the set of filling measured geodesic laminations
coincide with the Lebesgue measure up to scale.
Conformal densities induce ${\cal M}(S)$-invariant
Radon measures on measured lamination space (see the
discussion in Section 5). To understand ${\cal M}(S)$-invariant
Radon measures on ${\cal M\cal L}$ which
are \emph{not} of this form
we use train tracks as the main technical tool.
In this section we summarize
those properties of train tracks which are needed
for our purpose.

A \emph{maximal generic train track} on the surface $S$ is an embedded
1-complex $\tau\subset S$ whose edges
(called \emph{branches}) are smooth arcs with
well-defined tangent vectors at the endpoints. At any vertex
(called a \emph{switch}) the incident edges are mutually tangent.
Every switch is trivalent.
Through each switch there is a path of class $C^1$
which is embedded
in $\tau$ and contains the switch in its interior. In
particular, the branches which are incident
on a fixed switch are divided into
``incoming'' and ``outgoing'' branches according to their inward
pointing tangent at the switch.
The complementary regions of the
train track are trigons, i.e. discs with three cusps
at the boundary, or once punctured monogons,
i.e. once punctured discs with one cusp at the boundary.
We always identify train
tracks which are isotopic (see \cite{PH92} for a comprehensive
account on train tracks).

A maximal generic
train track or a geodesic lamination $\sigma$
is \emph{carried} by a train track $\tau$ if there is a
map $F:S\to S$ of class $C^1$ which is homotopic to the identity
and maps $\sigma$ into $\tau$ in such a way that the restriction
of the differential of $F$ to the tangent space of $\sigma$
vanishes nowhere; note that this makes sense since a train track
has a tangent line everywhere. We call the restriction of $F$ to
$\sigma$ a \emph{carrying map} for $\sigma$. Write $\sigma\prec
\tau$ if the train track or the geodesic lamination $\sigma$
is carried by the train track
$\tau$.

A \emph{transverse measure} on a maximal generic
train track $\tau$ is a
nonnegative weight function $\mu$ on the branches of $\tau$
satisfying the \emph{switch condition}:
For every switch $s$ of $\tau$, the sum of the weights
over all incoming branches at $s$
is required to coincide with the sum of
the weights over all outgoing branches at $s$.
The train track is called
\emph{recurrent} if it admits a transverse measure which is
positive on every branch. We call such a transverse measure $\mu$
\emph{positive}, and we write $\mu>0$.
The space ${\cal V}(\tau)$ of all transverse measures
on $\tau$ has the structure of a convex euclidean cone.
Via a carrying map, a measured geodesic lamination
carried by $\tau$ defines a transverse measure on
$\tau$, and every transverse measure arises in this
way \cite{PH92}. Thus ${\cal V}(\tau)$ can
naturally be identified with a
subset of ${\cal M\cal L}$ which is invariant under scaling.
A maximal generic train track
$\tau$ is recurrent if and only if the
subset ${\cal V}(\tau)$ of ${\cal M\cal L}$ has nonempty
interior.

A \emph{tangential measure} $\mu$ for
a maximal generic train track $\tau$ associates to every
branch $b$ of $\tau$ a nonnegative weight
$\mu(b)$ such that for every complementary triangle with
sides $c_1,c_2,c_3 $ we have $\mu(c_i)\leq \mu(c_{i+1})+
\mu(c_{i+2})$ (indices are taken
modulo three). The space ${\cal V}^*(\tau)$ of
all tangential measures on $\tau$ has the structure
of a convex euclidean cone.
The maximal generic train track $\tau$ is called
\emph{transversely
recurrent} if it admits a tangential
measure $\mu$ which is positive on every branch \cite{PH92}.
There is a one-to-one correspondence between
the space of tangential measures on $\tau$ and
the space of measured geodesic laminations
which \emph{hit $\tau$ efficiently}
(we refer to \cite{PH92} for an explanation of this
terminology). With this
identification, the pairing ${\cal V}(\tau)\times
{\cal V}^*(\tau)\to [0,\infty)$ defined
by $(\mu,\nu)\to \sum_b\mu(b)\nu(b)$ is just the
intersection form on ${\cal M\cal L}$ \cite{PH92}.
A maximal generic
train track $\tau$ is called
\emph{complete} if it is recurrent and transversely recurrent.
In the sequel we identify train tracks
which are isotopic, and we denote by ${\cal T\cal T}$
the set of isotopy classes of all complete train tracks
on $S$.

A half-branch $\hat b$ in a complete train track $\tau$ incident on
a switch $v$ of $\tau$ is called
\emph{large} if every embedded arc of class $C^1$
containing $v$ in its interior
passes through $\hat b$. A half-branch which is not large
is called \emph{small}.
A branch
$b$ in a complete train track
$\tau$ is called
\emph{large} if each of its two half-branches is
large; in this case $b$ is necessarily incident on two distinct
switches, and it is large at both of them. A branch is called
\emph{small} if each of its two half-branches is small. A branch
is called \emph{mixed} if one of its half-branches is large and
the other half-branch is small (for all this, see \cite{PH92} p.118).

There is a simple way to modify a complete train track $\tau$
to another complete train track.
Namely, if $e$ is a large branch of $\tau$ then we can perform a
right or left \emph{split} of $\tau$ at $e$ as shown in 
the figure below.
Note that a right split at $e$ is uniquely
determined by the orientation of $S$ and does not
depend on the orientation of $e$.
Using the labels in the figure, in the case of a right
split we call the branches $a$ and $c$ \emph{winners} of the
split, and the branches $b,d$ are \emph{losers} of the split. If
we perform a left split, then the branches $b,d$ are winners of
the split, and the branches $a,c$ are losers of the split.

\begin{figure}[ht]
\includegraphics{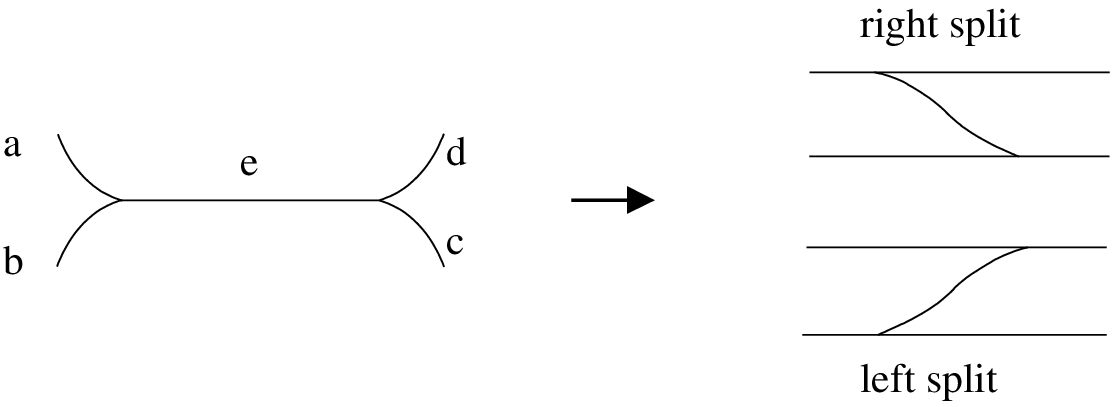}
\end{figure}

The split $\tau^\prime$ of a train track $\tau$ is carried
by $\tau$, and there is a natural choice of a carrying map which
maps the switches of $\tau^\prime$ to the switches of $\tau$.
There is a natural bijection of the set of branches
of $\tau$ onto the set of branches of $\tau^\prime$ which
maps the branch $e$ to the diagonal $e^\prime$ of the split.
The split of a complete train track is maximal,
transversely recurrent and generic, but it may not be
recurrent \cite{PH92}.

For a complete train track $\tau\in {\cal T\cal T}$ denote by
${\cal V}_0(\tau)\subset {\cal V}(\tau)$ the convex set of all
transverse measures on $\tau$ whose total weight (i.e. the sum of
the weights over all branches of $\tau)$ equals one. Let moreover
${\cal Q}(\tau)\subset {\cal Q}^1(S)$ be the set of all area one
quadratic differentials whose horizontal measured geodesic
lamination is contained in ${\cal V}_0(\tau)$ and whose vertical measured
geodesic lamination hits $\tau$ efficiently. Then ${\cal Q}(\tau)$
is a closed subset of ${\cal Q}^1(S)$ which however is not
compact.

Recall from \cite{MM99} that a \emph{vertex cycle} for $\tau$ is a
transverse measure on $\tau$ which spans an extreme ray in
the convex euclidean cone ${\cal
V}(\tau)$. Up to rescaling, such a vertex cycle is the counting
measure of a simple closed curve which is carried by $\tau$
(see \cite{MM99} for this fact). In
the sequel we identify a vertex cycle for $\tau$
with this simple closed
curve on $S$. With this convention, the transverse
measure on $\tau$ defined by a vertex cycle is integral.
Define a map $\Psi:{\cal T\cal T}\to {\cal C}(S)$ by
associating to a complete train track $\tau$ one of its vertex
cycles. Recall from Section 2 the definition
of the map $\Upsilon_{\cal Q}:{\cal Q}^1(S)\to
{\cal C}(S)$. We have.

\begin{lemma}\label{upsilon} There is a constant
$\chi_2>0$ such that
$d(\Upsilon_{\cal Q}(q),\Psi(\tau))\leq \chi_2$
for every $\tau\in {\cal T\cal T}$ and every
$q\in {\cal Q}(\tau)$.
\end{lemma}

\begin{proof} Let $\tau\in {\cal T\cal T}$ and let $q\in {\cal
Q}(\tau)$. Then the transverse measure
$\lambda\in {\cal V}_0(\tau)$ of the
horizontal measured geodesic
lamination  of $q$ can be written in
the form $\lambda=\sum_i a_i\alpha_i$ where $\alpha_i$ is the
transverse measure on $\tau$ defined by a vertex cycle and where
$a_i\geq 0$. Now the number of vertex cycles for $\tau$ is bounded
from above by a universal constant only depending on the
topological type of the surface $S$. The weight disposed on a branch of
$\tau$ by the counting measure of a vertex cycle does not exceed
two \cite{H06a}. This means that there is a number $a>0$ only
depending on the topological type of $S$ and there is some $i>0$
such that $a_i\geq a$.

The vertical measured geodesic lamination
of $q$ defines a tangential measure $\xi$ on $\tau$
with $\sum_b \lambda(b)\xi(b)=i(\lambda,\xi)=1$.
This implies that the vertex cycle $\alpha=\alpha_i$ satisfies
$i(\alpha,\xi)=\sum_{b}\alpha(b)\xi(b)\leq
i(\lambda,\xi)/a=1/a$. On the other hand,
the intersection number $i(\lambda,\alpha)$ is bounded
from above by a universal constant as well (Corollary
2.3 of \cite{H06a}).
Moreover, the $q$-length of the simple closed
curve $\alpha$ is bounded from above by
$2i(\alpha,\lambda)+2i(\alpha,\xi)$ (compare e.g. \cite{R05}
for this simple fact). Hence the $q$-length of $\alpha$
is uniformly bounded and therefore
Lemma \ref{simple1} implies that the distance between $\Upsilon_{\cal Q}(q)$
and $\alpha$ is bounded from above by a universal
constant. Since on the other hand the distance
in ${\cal C}(S)$ between any two
vertex cycles of a complete train track on $S$ is
also uniformly bounded \cite{H06a}, the lemma follows.
\end{proof}

Recall from \cite{MM99,H06b} the definition of a
\emph{framing} (or \emph{marking})
for the surface $S$. Such a framing $F$ consists of a pants
decomposition $P$ for $S$ and a collection of $3g-3+m$
essential simple
closed \emph{spanning curves}. For each pants curve $c\in P$
there is a unique such spanning curve which is disjoint from
$P-c$, which is not freely homotopic to any pants curve
of $P$ and
which has the minimal number
of intersection points with $c$ among all simple closed
curves with these properties. Any two such curves differ
by a multiple Dehn twist about $c$.
There is a number $\hat \chi_0>0$ and a hyperbolic
metric $h\in {\cal T}(S)$ such that the
$h$-length of each curve from our framing $F$
is at most $\hat \chi_0$. Such a metric can be constructed
as follows. Equip each pair of pants defined by $P$ with
a hyperbolic metric such that the length of each boundary
geodesic equals one. Then glue these pair of pants
in such a way that the spanning curves have the smallest
possible length.
We call such a hyperbolic metric
\emph{short for $F$}. By standard hyperbolic trigonometry,
there is a number $\epsilon>0$ such that
every hyperbolic metric which is short for some framing
$F$ of $S$ is contained in the set ${\cal T}(S)_\epsilon$
of all metrics whose \emph{systole},
i.e. the shortest length of a closed geodesic,
is at least $\epsilon$. Moreover,
the diameter in ${\cal T}(S)$
of the set of all hyperbolic metrics which are short
for a fixed framing $F$ is bounded from
above by a universal constant.

The mapping class group ${\cal M}(S)$
naturally acts on the collection of all framings of $S$.
By equivariance and the
fact that the number of orbits for the
action of ${\cal M}(S)$ on ${\cal T\cal T}$
is \emph{finite}, there is a number $k>0$
and for every complete
train track $\tau\in {\cal T\cal T}$ there is a framing
$F$ of $S$ which consists of simple closed curves carried
by $\tau$ and such that the
total weight of the counting measures on $\tau$
defined by these curves does not
exceed $k$ (compare the discussion in \cite{MM99}).
We call such a framing \emph{short for $\tau$}.
The intersection number $i(c,c^\prime)$ between any two
simple closed curves
$c,c^\prime$ which are carried by some
$\tau\in {\cal T\cal T}$ and which define counting measures
on $\tau$ of total weight at most $k$ is bounded from
above by a universal constant
(see Corollary 2.3 of \cite{H06a}).
In particular,
for any two short framings $F,F^\prime$
for $\tau$ the distance in ${\cal T}(S)$ between
any two hyperbolic metrics which are short for $F,F^\prime$
is uniformly bounded
\cite{MM99,Bw06}.

Define a map $\Lambda:{\cal T\cal T}\to
{\cal T}(S)$ by associating to a complete train track
$\tau$ a hyperbolic metric
$\Lambda(\tau)\in {\cal T}(S)$ which is short for a
short framing for $\tau$. By our above discussion, there
is a number $\chi_3>0$ only depending on the topological
type of $S$ such that if $\Lambda^\prime$ is another
choice of such a map then we have
$d(\Lambda(\tau),\Lambda^\prime(\tau))\leq \chi_3$ for every
$\tau\in {\cal T\cal T}$.
In particular, the map $\Lambda$ is \emph{coarsely
${\cal M}(S)$-equivariant}: For every $\tau\in {\cal T\cal T}$
and every $g\in {\cal M}(S)$ we have
$d(\Lambda(g\tau),g\Lambda(\tau))\leq \chi_3$.
The next observation will be useful in Section 5 (compare
also \cite{H05} for a related result). For its formulation,
let again $P:{\cal Q}^1(S)\to {\cal T}(S)$ be
the canonical projection.

\begin{lemma}\label{vertexcycle} There is a number $\ell >0$ and
for every $\epsilon >0$ there is a number
$m(\epsilon)>0$ with the following property.
Let $\sigma,\tau\in {\cal T\cal T}$ and
assume that $\sigma$ is carried by $\tau$ and that the distance
in ${\cal C}(S)$ between any vertex cycle of $\sigma$
and any vertex cycle of $\tau$ is at least $\ell$.
Let $q\in {\cal Q}(\tau)$ be a quadratic
differential whose
horizontal measured geodesic lamination $q_h$ is carried by
$\sigma$.
If the total weight of the transverse
measure on $\sigma$ defined by $q_h$
is not smaller than $\epsilon$, then
$d(\Lambda(\tau),Pq)\leq m(\epsilon)$.
\end{lemma}
\begin{proof}
Let $\chi_0>0$ be as in the definition of the
map $\Upsilon_{\cal T}:{\cal T}(S)\to {\cal C}(S)$ in Section 2.
Let $\chi_1>0$ be as in Lemma \ref{simple2} and let
$\kappa_1>0$ be
such that for every
$x\in {\cal T}(S)$ the diameter in ${\cal C}(S)$ of
the set of simple closed curves whose $x$-length is at
most $\chi_0$ is bounded from above by $\kappa_1$.
Let $\chi_2>0$ be
as in Lemma \ref{upsilon} and let $\ell >2\kappa_1+2\chi_1+2\chi_2+3$.
Let $\sigma\prec\tau\in {\cal T\cal T}$
be such that the distance in ${\cal C}(S)$ between
any vertex cycle of $\tau$ and any vertex cycle
of $\sigma$ is at least $\ell$. Let $\epsilon >0$ and let
$q\in {\cal Q}(\tau)$ be such that the horizontal
measured geodesic lamination $q_h$ of $q$ is carried
by $\sigma$ and is such that the total weight
of the transverse measure on $\sigma$ defined by
$q_h$ is not smaller than $\epsilon$.
We claim that the $Pq$-length of any essential simple
closed curve on $S$ is at least $\epsilon\chi_0$.
Namely, if $t\geq 0$ is such that $\Phi^tq\in {\cal Q}(\sigma)$
then $t\leq -\log\epsilon$. Thus if there is a simple
closed curve $c$ on $S$ whose $Pq$-length is smaller than
$\epsilon\chi_0$ then the $P\Phi^tq$-length of
$c$ is at most $\chi_0$ (see \cite{IT99}
for this well known result of Wolpert).
By the choice of the constant
$\kappa_1$, this means that $d(\Upsilon_{\cal T}(Pq),
\Upsilon_{\cal T}(P\Phi^tq))\leq 2\kappa_1$ and hence
from Lemma \ref{simple2} we conclude that
$d(\Upsilon_{\cal Q}(q),
\Upsilon_{\cal Q}(\Phi^tq)\leq 2\kappa_1+2\chi_1$.
From Lemma \ref{upsilon} we conclude that
$d(\Psi(\sigma),\Psi(\tau))\leq 2\kappa_1+2\chi_1+2\chi_2$
which contradicts our choice of $\ell$.
In particular, we have $Pq\in
{\cal T}(S)_{\epsilon\chi_0}$.

Now by Lemma 3.3 of \cite{Mi94}, there is a constant
$L>0$ depending on $\epsilon$ such that for
every $q\in {\cal Q}^1(S)$ with
$Pq\in {\cal T}(S)_{\epsilon\chi_0}$
the singular euclidean metric defined by
$q$ is $L$-bilipschitz equivalent to the
hyperbolic metric $Pq$. Thus for the
proof of our lemma it suffices to show
that
for a quadratic differential $q\in {\cal Q}(\tau)$ as
above the $q$-length of any short framing for $\tau$
is uniformly bounded.

Let $q_v$ be the vertical measured geodesic
lamination of $q$. Then $q_v$ defines a tangential
measure on $\tau$. We claim that the
weight of $q_v$ disposed on any branch of $\tau$ is uniformly
bounded. For this note that the
transverse measure on the train track
$\sigma$ defined by the horizontal measured
geodesic lamination $q_h$ of $q$ can be reprensented in the form
$q_h=\sum_i a_i\alpha_i$ with vertex cycles
$\alpha_i$ for $\sigma$. Now a counting measure of a simple closed
curve which is carried by $\sigma$ is integral, and the total weight
of the counting measure defined by a vertex cycle of $\sigma$ is
bounded from above by a universal constant
$k>0$. Moreover, the
number of vertex cycles for $\sigma$ is bounded from above by a
universal constant $p>0$. Since by assumption the total
weight of $q_h$ on $\sigma$ is at least $\epsilon$, we have
$a_i\geq \epsilon/kp$ for at least one $i$.

By the choice of the constant $\ell>3$, the distance in ${\cal
C}(S)$ between any vertex cycle of $\sigma$ and any
vertex cycle for $\tau$ is at least $\ell$.
Thus by Lemma 4.9 of \cite{MM99}, the image of $\alpha_i$
under a carrying map $\alpha_i\to \tau$ is a \emph{large
subtrack} $\omega$ of $\tau$. This means that this image is
a train track on $S$ which is a subset of $\tau$ and whose
complementary components are topological discs or
once punctured topological discs. Thus each
such complementary
component is a polygon or a once punctured polygon which is
a union of complementary components of $\tau$.

The map ${\cal V}(\sigma)\to {\cal V}(\tau)$
is convex linear and therefore the
$q_h$-weight of every branch of $\omega$ is bounded from below by
$\epsilon/kp$. In particular, if $q_v$ denotes the vertical
measured geodesic lamination of $q$ then we obtain from the
identity $i(q_h,q_v)=1$ that the weight of every branch of $\omega$
with respect to the tangential measured defined by $q_v$ is
bounded from above by $kp/\epsilon$. From the definition of a
tangential measure for $\tau$ and the fact that $\omega$ is large
we deduce that the total weight that $q_v$ disposes on the branches
of $\tau$ is bounded from above by a constant only
depending on $\epsilon$. Namely, let $D$ be a complementary
polygon of $\omega$ with more than 3 sides. Then
the branches of $\tau$ which are contained in $D$
decompose $D$ into triangles. There is at least one such triangle
$T$ with two sides contained in the boundary of $D$,
i.e. with two sides contained in $\omega$.
By the definition of a tangential measure,
the total weight disposed by $q_v$
on the third side of $T$ which is contained
in the interior of $D$ is bounded from above by the sum of the
total weights on the sides of $T$ contained
in $\omega$. Thus by the
above consideration, the total weight disposed by $q_v$ on
the boundary of $T$
is uniformly bounded. Since the number of complementary
components of $\tau$ only depends on the topological type
of $S$, with a uniformly bounded number of steps we
obtain inductively the above claim (note that the argument
is also valid for once punctured complementary polygons of $\omega$).

Since the total weight of the tangential measure
on $\tau$ defined by $q_v$ is uniformly bounded,
the intersection number $i(c,q_v)$ for every
curve from a short framing for $\tau$ is uniformly
bounded as well. Moreover,
by Corollary 2.3 of \cite{H06a}
there is a constant $\kappa_2>0$ with the
following property. Let $\tau\in {\cal T\cal T}$
and let $c$ be a curve contained in a short framing for
$\tau$. Then for every $\nu\in {\cal V}_0(\tau)$
we have $i(\nu,c)\leq \kappa_2$.
Now for every quadratic differential $z\in {\cal Q}^1(S)$
with horizontal and vertical measured
geodesic lamination $z_h,z_v$ the $z$-length
of a simple closed curve $c$
is bounded from above by $2i(z_h,c)+2i(z_v,c)$.
This shows that the $q$-length of every simple closed curve
from a short framing of $\tau$ is uniformly bounded and
completes the proof of the lemma.
\end{proof}

Every projective measured geodesic lamination
$\xi\in {\cal P\cal M\cal L}$ determines
an unstable manifold $W^u(\xi)\subset {\cal Q}^1(S)$
of all area one quadratic differentials
whose vertical measured geodesic lamination
is contained in the class $\xi$. This unstable manifold
can naturally be identified with the set
of all measured geodesic laminations $\zeta$
such that $\zeta$ and $\xi$ jointly
fill up $S$. By the Hubbard-Masur
theorem, the unstable manifold $W^u(\xi)$ projects
homeomorphically onto ${\cal T}(S)$ and hence
the Teichm\"uller metric defines a distance function
$d$ on $W^u(\xi)$.

For a train track $\tau\in {\cal T\cal T}$
let ${\cal P\cal E}(\tau)$ be the set of all
projective measured geodesic laminations which
hit $\tau$ efficiently. Note that if $\sigma\prec\tau$
then ${\cal P\cal E}(\tau)\subset {\cal P\cal E}(\sigma)$.
For
$\xi\in {\cal P\cal E}(\tau)$ and a number
$R>0$ the train track $\tau$ is called
\emph{$R-\xi$-tight} if the diameter
of ${\cal Q}(\tau)\cap W^u(\xi)$ with respect to the
lift of the Teichm\"uller metric is at most $R$.
The next corollary will be used in Section 5.

\begin{corollary}\label{tight}
Let $\ell >0$ be as in Lemma \ref{vertexcycle}.
Then for every $\epsilon >0$ there is a number
$R=R(\epsilon)>0$ with the following property.
Let $\eta\prec\sigma\prec\tau$ and assume that
the distance in ${\cal C}(S)$ between any
vertex cycle of $\sigma$ and any vertex cycle of
$\eta$ as well as any vertex cycle of $\tau$ is at least $\ell$.
Let $q\in {\cal Q}(\tau)$ be a quadratic differential
whose horizontal measured geodesic lamination $q_h$
is carried by $\eta$. Let $\xi\in {\cal P\cal E}(\tau)$
be the projective class of the vertical measured
geodesic lamination of $q$. If the total weight
of the transverse measure on $\eta$ defined
by $q_h$ is not smaller than $\epsilon$ then
$\sigma$ is $R-\xi$-tight.
\end{corollary}
\begin{proof}
Let $\delta>0$ be such that
$\Lambda(\tau)\in {\cal T}(S)_\delta$ for every
$\tau\in {\cal T\cal T}$. By possibly decreasing
$\delta$ we may assume that ${\cal T}(S)_\delta$
is connected.
If we equip
${\cal T}(S)_\delta$ with the length metric
induced by the Finsler structure defining the
Teichm\"uller metric then ${\cal T}(S)_\delta$ is
a proper geodesic metric space on which the mapping class
group ${\cal M}(S)$ acts
properly discontinuously and cocompactly as a group
of isometries.
The set ${\cal T\cal T}$ is the set of vertices
of the \emph{train track complex} \cite{H06b} which is
a connected metric graph on which the mapping class
group acts properly and cocompactly as a group of
isometries. Since the map $\Lambda$ is coarsely
${\cal M}(S)$-equivariant this means that there is some
$L>1$ such that
$\Lambda:{\cal T\cal T}\to {\cal T}(S)_\delta$
is an $L$-quasi-isometry.

Let $\eta\prec\sigma\prec\tau$ be
as in the corollary, let $\epsilon >0$ and
assume that there is some $q\in {\cal Q}(\tau)$ such that
the horizontal measured geodesic lamination
$q_h$ of $q$ is carried by $\eta$ and that the weight
disposed on $\eta$ by $q_h$ is bounded from below by
$\epsilon$. Then the weight $c$ disposed on $\sigma$
by $q_h$ is contained in the interval $[\epsilon,1]$.
By Lemma \ref{vertexcycle}, applied both
to the train tracks $\sigma\prec \tau$ and
to the train tracks $\eta\prec\sigma$ (with the
quadratic differential $\Phi^sq$ for $s=-\log c$)
the distance between $\Lambda(\sigma)$ and $\Lambda(\tau)$
is bounded from above by a number $\rho >0$ only
depending on $\epsilon$. Then the distance
between $\Lambda(\sigma)$ and $\Lambda(\tau)$ in
${\cal T}(S)_\delta$ is uniformly bounded as well
and hence the same is true for the distance between
$\sigma,\tau$ in the train track complex. However
there are only finitely many orbits under the
action of the mapping class group of pairs
$\sigma\prec\tau$ whose distance in the train track
complex is uniformly bounded. Thus
by invariance under the mapping class group,
there is a universal number $\kappa >0$ with the following
property. If $\alpha$ is any vertex cycle of
$\sigma$, then the total weight disposed on
$\tau$ by $\alpha$ is at most $\kappa$.
Thus the carrying map $\sigma\to \tau$
maps ${\cal V}_0(\sigma)$ to a subset of
${\cal V}(\tau)$ consisting of transverse measures
whose total weight is bounded from above by a universal constant.

As a consequence, there is a universal constant
$\rho >0$ such that
if $z\in {\cal Q}(\tau)\cap W^u(\xi)$ is
\emph{any} quadratic differential with the property
that the horizontal measured geodesic lamination $z_h$ of
$z$ is carried by $\sigma$ then the total
weight which is disposed by $z_h$ on $\sigma$ is
bounded from below by $\rho$.
By Lemma \ref{vertexcycle}, the
the set of all quadratic differentials
$q\in {\cal Q}(\tau)\cap W^u(\xi)$
whose horizontal measured geodesic lamination is carried by
$\sigma$ projects to a ball of uniformly bounded radius
about $\Lambda(\tau)$. Together with the above observations
this shows the corollary.
\end{proof}

\section{Invariant Radon measures on ${\cal M\cal L}$}

In this section we complete the proof of the theorem from
the introduction. We continue to use the
assumptions and notations from Section 3.
Recall first that for every point
$x\in {\cal T}(S)$ and every $\xi\in {\cal P\cal M\cal L}$
there is a unique quadratic differential
$q(x,\xi)\in {\cal Q}^1(S)_x$ of area one
on the Riemann surface
$x$ whose horizontal measured geodesic lamination
$q_h(x,\xi)$ is contained in the class $\xi$.
The assignment
$\xi\to q_h(x,\xi)$ determines a
homeomorphism
${\cal P\cal M\cal L}\times
\mathbb{R}\to {\cal M\cal L}$ by assigning 
to $(\xi,t)\in
{\cal P\cal M\cal L}\times \mathbb{R}$ the
measured geodesic lamination $e^tq_h(x,\xi)$.

A conformal density $\{\nu^y\}$ on ${\cal P\cal M\cal L}$
of dimension $\alpha$
defines a Radon measure $\Theta_\nu$ on
${\cal M\cal L}$ via
$d\Theta_\nu(\xi,t)=d\nu^x(\xi)\times e^{\alpha t}dt$
where $\xi\in {\cal P\cal M\cal L}, t\in \mathbb{R}$.
By construction, this measure
is quasi-invariant under the one-parameter group of
translations $T^s$ on ${\cal M\cal L}={\cal P\cal M\cal L}\times
\mathbb{R}$
given by $T^s(\xi,t)=(\xi,s+t)$. More precisely, we have
$\frac{d\nu\circ T^s}{d\nu}=e^{\alpha s}$.

The measure $\Theta_\nu$ is moreover invariant under the action of
the mapping class group ${\cal M}(S)$.
Namely, for $\xi\in {\cal P\cal M\cal L}$
and $g\in {\cal M}(S)$
the lamination
$g(q_h(x,\xi))=q_h(g(x),g(\xi))$
equals $e^{\Psi(x,g(x),g(\xi))}
q_h(x,g(\xi))$ where $\Psi$ is the cocycle defined
in the beginning of Section 3.
On the other hand, we have
\begin{equation}
d (\Theta_\nu\circ g)(\xi,t)
=d\nu^{g(x)}(g(\xi))
\times e^{\alpha t}dt.\end{equation}
By the
definition of a conformal density, for
$\nu^x$-almost every $\xi\in {\cal P\cal M\cal L}$
the Radon Nikodym derivative of the measure
$\nu^{g(x)}$ with respect to $\nu^x$ at the point
$g(\xi)$
equals $e^{\alpha \Psi(x,g(x),g(\xi))}$ and therefore the
measure $\Theta_\nu$ is indeed invariant
under the action of ${\cal M}(S)$.
As a consequence, every conformal density on
${\cal P\cal M\cal L}$ induces a ${\cal M}(S)$-invariant
Radon measure on ${\cal M\cal L}={\cal P\cal M\cal L}\times
\mathbb{R}$ which is quasi-invariant under the one-parameter
group of translations $T^s$.

Now let $\eta$ be any ergodic ${\cal M}(S)$-invariant
Radon measure on ${\cal M\cal L}$.
Then
\begin{equation}
H_\eta=\{a\in \mathbb{R}\mid \eta\circ T^a\sim \eta\}
\end{equation}
is a closed subgroup of $\mathbb{R}$ \cite{ANSS02}.
The next lemma is an easy consequence of Proposition
\ref{confonfilling} and
Sarig's cocycle reduction theorem \cite{S04}.
For its formulation, recall from Section 3 the
definition of the space ${\cal F\cal M\cal L}$
of all projective measured geodesic laminations
which fill up $S$. We call a measured geodesic
lamination $\lambda\in {\cal M\cal L}$
\emph{filling} if its projectivization
is contained in ${\cal F\cal M\cal L}$.
The set of all filling measured geodesic laminations
is invariant under the mapping class group.

\begin{lemma}\label{groupnontrivial}
Let $\eta$ be an ergodic ${\cal M}(S)$-invariant
Radon measure on ${\cal M\cal L}$ which
gives full mass to the filling measured geodesic
laminations. If $H_\eta\not=\{0\}$ then $\eta$
coincides with the Lebesgue measure up to scale.
\end{lemma}
\begin{proof}
Define a measurable countable
equivalence relation ${\cal R}$ on
${\cal P\cal M\cal L}$ by
$\chi {\cal R} \xi$ if and only if
$\chi$ and $\xi$ are contained in the
same orbit for the action of the mapping
class group. Recall the definition
of the cocycle
$\Psi:{\cal T}(S)\times{\cal T}(S)\times
{\cal P\cal M\cal L}\to \mathbb{R}$.
For a fixed point $x\in {\cal T}(S)$ we obtain
a real-valued cocycle for the
action of ${\cal M}(S)$ on ${\cal P\cal M\cal L}$,
again denoted by $\Psi$,
via $\Psi(\lambda,g)=\Psi(x,g^{-1}x,\lambda)$
$(\lambda\in
{\cal P\cal M\cal L}$ and $g\in {\cal M}(S)$).
By the cocycle identity (\ref{cocycle}) for $\Psi$ we have
$\Psi(\lambda,hg)=\Psi(\lambda,g)+
\Psi(g\lambda,h)$, i.e.
$\Psi$ is indeed a cocycle which
can be viewed as a cocycle on ${\cal R}$.
We also write $\Psi(\lambda,\xi)$ instead
of $\Psi(\lambda,g)$
whenever $\xi=g\lambda$; note that
this is only well defined if $\lambda$ is
not fixed by any element of ${\cal M}(S)$, however
this ambiguity will be of no importance in the
sequel.

Recall that the choice of a point $x\in {\cal T}(S)$
determines a homeomorphism ${\cal M\cal L}\to
{\cal P\cal M\cal L}\times \mathbb{R}$.
The cocycle $\Psi$ then
defines an equivalence relation
${\cal R}_\Psi$ on ${\cal M\cal L}
={\cal P\cal M\cal L}\times \mathbb{R}$
by
\begin{equation}{\cal R}_\Psi=\{((\lambda,t),(\xi,s))\in
({\cal P\cal M\cal L}\times \mathbb{R})^2\mid
(\lambda,\xi)\in {\cal R}\text{ and } s-t=
\Psi(\lambda,\xi)\}.\end{equation}

Let $\eta$ be an ergodic ${\cal M}(S)$-invariant
Radon measure on ${\cal M\cal L}$ which gives full
mass to the filling measured geodesic laminations.
By the results in \cite{ANSS02},
if $H_\eta=\mathbb{R}$ then $\eta$ is induced
by a conformal density as in
the beginning of this section. In particular,
by Proposition \ref{confonfilling},
in this case the
measure $\eta$ equals the Lebesgue measure up to scale.
Thus for the proof of
our lemma we are left with the
case that
$H_\eta=c\mathbb{Z}$ for a number
$c> 0$.

By the cocycle reduction theorem of Sarig
(Theorem 2 of \cite{S04}), in this case there is a Borel
function $u:{\cal P\cal M\cal L}\to \mathbb{R}$
such that $\Psi_u(x,y)=\Psi(x,y)+u(y)-u(x)\in H_\eta$
holds $\eta$-almost everywhere in ${\cal R}_\Psi$.
Since $c>0$ we may assume
without loss of generality that the
function $u$ is \emph{bounded}.
Following \cite{S04}, for $a\in \mathbb{R}$ define
$\theta_a(x,t)=(x,t-u(x)-a)$.
By Lemma 2 of \cite{S04},
for a suitable choice of $a$ the
measure $\eta\circ \theta_a^{-1}$ is an
${\cal R}_{\Psi_u}$-invariant ergodic Radon
measure supported on ${\cal P\cal M\cal L}\times c\mathbb{Z}$.

We now follow Ledrappier and Sarig \cite{LS06}
(see also \cite{Ld06}). Namely,
since $\eta$ is invariant and
ergodic under the action of ${\cal M}(S)$ and since
the $\mathbb{R}$-action
on ${\cal M\cal L}$ commutes with the
${\cal M}(S)$-action, for every $t\in \mathbb{R}$
the measure
$\eta\circ T^t$ is also ${\cal M}(S)$-invariant
and ergodic.
Thus either $\eta\circ T^t$ and $\eta$ are
singular or they coincide up to scale.
As a consequence, there is some number
$\alpha \in \mathbb{R}$ such that $\eta\circ T^c=
e^{\alpha c}\eta$. Since $\theta=\theta_a$ and $T^t$ commute,
we also have $\eta\circ (\theta^{-1}\circ T^c)=
e^{\alpha c}\eta\circ \theta^{-1}$.
Consequently the measure $e^{-\alpha t}\eta\circ\theta^{-1}$
is invariant under the translation $T^c$.
Since moreover
$e^{-\alpha t}\eta\circ \theta^{-1}$ is supported
in ${\cal P\cal M\cal L}\times c\mathbb{Z}$, it follows
that $e^{-\alpha t}\eta\circ \theta^{-1}=
\nu\times m_{H_\eta}$ with some measure
$\nu$ on ${\cal P\cal M\cal L}$.

The measure $\nu$ is
finite since $\eta\circ\theta^{-1}$ is Radon and the function
$u$ is bounded.
The measure $\eta$ is ${\cal M}(S)$-invariant
and therefore $\eta\circ \theta^{-1}$
is invariant under $\theta\circ{\cal M}(S)\circ
\theta^{-1}$. In particular, the
measure class of $\nu$ is invariant under
the action of ${\cal M}(S)$. More precisely, we have
\begin{equation}
\frac{d\nu\circ g}{d\nu}(\xi)=e^{\alpha\Psi(\xi,g)}
\frac{e^{-\alpha u(\xi)}}{e^{-\alpha u(g(\xi))}}\end{equation}
for all $g\in {\cal M}(S)$ and
$\nu$-almost every $\xi\in {\cal P\cal M\cal L}$
(see \cite{LS06}).
As a consequence, if we define
$d\nu^x(\xi)=e^{\alpha u(\xi)}d\nu(\xi)$ and
$d\nu^{y}(\xi)=e^{\alpha \Psi(x,y,\xi)}d\nu^x(\xi)$ for
$y\in {\cal T}(S)$ then $\{\nu^y\}$ defines
a conformal density of dimension $\alpha$ on
${\cal P\cal M\cal L}$. Note that the measure
$\nu^x=e^{\alpha u}\nu$
is finite since the function $u$ is bounded.
By our assumption on the measure
$\eta$, the conformal density $\{\nu^x\}$
gives full measure to the ${\cal M}(S)$-invariant set
${\cal F\cal M\cal L}$ of projective measured
laminations which fill up $S$
and hence we conclude from
Proposition \ref{confonfilling} that $\eta$ equals
the Lebesgue measure $\lambda$ up to scale.
However, the Lebesgue measure is quasi-invariant under
the translations $\{T^t\}$ which
is a contradiction
to the assumption that $H_\eta=c\mathbb{Z}$ for some
$c>0$. This shows the lemma.
\end{proof}

The investigation of
${\cal M}(S)$-invariant ergodic measures $\eta$ on ${\cal M\cal L}$
which give full measure to the filling laminations and
satisfy $H_\eta=\{0\}$ is more difficult.
We begin with an observation which is similar to
Proposition \ref{confonfilling}.
For this call 
a measured geodesic lamination \emph{weakly recurrent}
if its projectivization is contained in the set 
${\cal R\cal M\cal L}$.

\begin{lemma}\label{fillingrec}
An ${\cal M}(S)$-invariant Radon measure $\eta$ on ${\cal M\cal L}$
which gives full measure to the filling
measured geodesic laminations gives full measure
to the recurrent measured geodesic laminations.
\end{lemma}
\begin{proof}
Let $\eta$ be an ${\cal M}(S)$-invariant ergodic
Radon measure on ${\cal M\cal L}$
which gives full measure to the filling measured geodesic laminations.
We use the measure $\eta$ to construct a
locally finite Borel
measure $\nu$ on ${\cal Q}(S)$ which is invariant
under the horocycle flow. Namely,
for every $q\in {\cal Q}^1(S)$ the assignment which
associates to a quadratic differential $z$ contained
in the unstable manifold $W^u(q)$ its horizontal
measured geodesic lamination is a homeomorphism of
$W^u(q)$ into ${\cal M\cal L}$. Thus by equivariance
under the action of the mapping class group, the measure
$\eta$ lifts to a ${\cal M}(S)$-invariant
family $\{\eta^u\}$  of locally finite
measures on unstable manifolds $W^u(q)$
$(q\in {\cal Q}^1(S))$.
This family of measures then projects to a family
$\{\nu^u\}$ of locally finite Borel measures on
the leaves of the unstable foliation on ${\cal Q}(S)$.
The family $\{\nu^u\}$ is invariant under holonomy
along strong stable manifolds.

Define $d\nu=d\nu^u\times d\lambda^{ss}$ where
$\lambda^{ss}$ is a standard family of
Lebesgue measures on strong stable manifolds
which is invariant under the horocycle flow $h_t$.
Then $\nu$ is a locally
finite $h_t$-invariant Borel measure on ${\cal Q}(S)$.
For $\nu$-almost every point $q\in {\cal Q}(S)$ the
horizontal and the vertical measured geodesic
laminations of $q$ fill up $S$.
As in the proof of Proposition \ref{confonfilling},
Dani's argument together with the theorem from the appendix
implies that $\nu$ is \emph{finite}. 
Moreover, for every $\epsilon >0$ there is a compact
subset $K_\epsilon$ of ${\cal Q}(S)$ not depending on $\nu$ such that
$\nu(K_\epsilon)/\nu({\cal Q}(S))\geq 1-\epsilon$.

Via replacing $\nu$ by $\nu/\nu({\cal Q}(S)$ we
may assume that $\nu$ is a probability measure.
For $s>0$ define
\[\nu(s)=\frac{1}{s}\int_0^s\Phi^t\nu dt.\]
Since $h_t\circ \Phi^s=\Phi^s\circ h_{e^st}$ for all
$s,t\in \mathbb{R}$, the Borel probability measure
$\nu(s)$ is $h_t$-invariant and 
gives full measure to the points with filling horizontal
measured geodesic laminations. Therefore we have 
$\nu(s)(K_\epsilon)\geq 1-\epsilon$ for all $s>0$, all $\epsilon >0$.
This implies that 
there is a sequence $s_i\to \infty$ such that
the measures $\nu(s_i)$ converge as $i\to \infty$ weakly
to a Borel probability measure $\nu(\infty)$ on 
${\cal Q}(S)$ 
which is invariant
under both the horocycle flow and the Teichm\"uller geodesic
flow. By the Poincar\'e
recurrence theorem, $\nu(\infty)$ gives full measure to the
forward \emph{recurrent} quadratic 
differentials. 

As a consequence, $\nu$-almost every $q\in {\cal Q}(S)$
contains a forward recurrent point in
its $\omega$-limit set. Namely, for $\epsilon >0$
there is a compact subset $B_\epsilon$ of ${\cal Q}(S)$
which consists of forward recurrent points and such that
$\nu(\infty)(B_\epsilon)>1-\epsilon$. Let $\{U_\ell\}$ 
be a family of open
neighborhoods of $B_\epsilon$ such that $U_\ell\supset U_{\ell +1}$
for all $\ell$ and $\cap_\ell U_\ell=B_\epsilon$.
Then for every $\ell >0$ there is some $i(\ell)>0$ such that
$\nu(s_i)(U)\geq 1-2\epsilon$ for all $i\geq i(\ell)$.

For $\ell >0$ define $C_\ell=\{q\mid \Phi^{t}q\in U_\ell$
for infinitely many $t>0\}$; then $C_\ell\supset C_{\ell +1}$
for all $\ell$. We claim that 
$\nu(C_\ell)\geq 1-5\epsilon$ for every $\ell$. Namely,
otherwise there is a number $T>0$ and there is a subset
$A$ of ${\cal Q}(S)$ with $\nu(A)\geq 4\epsilon$ and such that
$\Phi^tz\not\in U_\ell$ for every $z\in A$ and every $t\geq T$. 
Then necessarily $\nu(s_i)(U_\ell)\leq 1-3\epsilon$ for all
sufficiently large $i$ which is impossible.
Since the neighborhood $U_\ell$ of $B_\epsilon$ was arbitrary
and since $\nu$ is Borel regular we conclude that
$\nu(\cap_\ell C_\ell)\geq 1-5\epsilon$. Thus the $\nu$-mass of 
all points $q\in {\cal Q}(S)$ which 
contain a point $z\in B_\epsilon$ in its $\omega$-limit set
is at least $1-5\epsilon$.
Since $B_\epsilon$ consists of recurrent points and
since $\epsilon >0$ was arbitrary we conclude that
$\nu$-almost every $q\in {\cal Q}(S)$ contains a forward
recurrent point in its $\omega$-limit set. 
This shows the lemma.
\end{proof}

For the analysis of ${\cal M}(S)$-invariant Radon measures
$\eta$ on ${\cal M\cal L}$ with $H_{\eta}=\{0\}$
we use a construction reminiscent of symbolic dynamics where
the Markov shift is replaced by complete train tracks and
their splits. We next establish some technical preparations
to achieve this goal.

Define a geodesic lamination $\xi$ on $S$ to
be \emph{complete} if $\xi$ is maximal
and can be approximated in the Hausdorff
topology by simple closed geodesics.
The space ${\cal C\cal L}$ of all complete
geodesic laminations equipped with the
Hausdorff topology is a compact ${\cal M}(S)$-space
\cite{H06b}.
A transversely recurrent generic train track is complete
if and only if it carries a complete geodesic
lamination. Every minimal geodesic lamination $\lambda$
is a sublamination of a complete geodesic lamination.
If $\lambda$ fills up $S$ then the number
of complete geodesic laminations which contain
$\lambda$ as a sublamination is bounded from
above by a universal constant (for all this,
see \cite{H06b}).

Denote by ${\cal H}\subset {\cal C\cal L}$ the set of all
complete geodesic laminations $\lambda\in {\cal C\cal L}$
which contain a uniquely ergodic minimal component
which fills up $S$.
The mapping class group ${\cal M}(S)$ naturally
acts on ${\cal H}$ as a group of transformations.
There is a
finite-to-one ${\cal M}(S)$-equivariant
map \begin{equation}
E:{\cal H}\to {\cal F\cal M\cal L}\end{equation}
which associates
to every $\lambda\in {\cal H}$ the unique projective
measured geodesic lamination which is supported in $\lambda$.
The number of preimages in ${\cal H}$ of a point in
${\cal F\cal M\cal L}$ is bounded from above by a universal constant.

A \emph{full split} of a
complete train track $\tau$ is a complete train track
$\sigma$ which can be obtained from $\tau$ by splitting
$\tau$ at each large branch precisely once.
A \emph{full splitting sequence} is a sequence
$\{\tau_i\}\subset {\cal T\cal T}$
such that for each $i$, the train track
$\tau_{i+1}$ can be obtained from $\tau_{i}$ by a full split.
For a complete train track $\tau$ denote by
${\cal C\cal L}(\tau)$ the set of all complete
geodesic laminations which are carried by $\tau$. Then
${\cal C\cal L}(\tau)$ is a subset of ${\cal C\cal L}$ which
is both open and closed.
If $\{\tau_i\}$ is an infinite full splitting sequence
then $\cap_i {\cal C\cal L}(\tau_i)$ consists of a unique
point. For every
complete train track
$\tau$ and every complete geodesic
lamination $\lambda\in {\cal C\cal L}(\tau)$ there is
a unique full splitting sequence
$\{\tau_i(\lambda)\}$ issuing
from $\tau_0(\lambda)=\tau$
such that $\cap_i {\cal C\cal L}(\tau_i(\lambda))=\{\lambda\}$
(for all this, see \cite{H06b}).

For $\tau\in {\cal T\cal T}$ let ${\cal
V}_0(\tau)\subset {\cal M\cal L}$ be the set of all measured
geodesic laminations $\nu$ whose support is carried by $\tau$ and
such that the total mass of the transverse measure on $\tau$
defined by $\nu$ equals one. If the complete geodesic lamination
$\lambda\in {\cal H}$ is carried by $\tau$ then there is a unique
measured geodesic lamination $\nu(\lambda,\tau)\in
{\cal V}_0(\tau)$
whose support is contained in $\lambda$.
There is a number
$a>0$ only depending on the topological type of $S$ such that for
all $\lambda\in {\cal H}\cap
{\cal C\cal L}(\tau)$ and all $i\geq 0$ we have
$\nu(\lambda,\tau_{i+1}(\lambda))=e^s\nu(\lambda,
\tau_i(\lambda))$ for some $s\in [0,a]$.
Namely, $\tau_{i+1}(\lambda)$ is obtained from
$\tau_i(\lambda)$ by a uniformly bounded number
of splits. Moreover,
if $\eta\in {\cal T\cal T}$ is obtained from
$\tau$ by a split at a large branch $e$, with losing
branches $b,d$, and if $\nu$ is any transverse measure
on $\eta$, then $\nu$ projects to a transverse measure
$\hat \nu$
on $\tau$ with the following properties. Using the natural
identification of the branches of $\tau$ with the
branches of $\eta$, for
every branch $h\not=e$ of $\tau$, the
$\hat\nu$-weight of $h$ coincides with the $\nu$-weight
of the branch $h$ in $\eta$.
Moreover, we have $\hat \nu(e)=\nu(e)+\nu(b)+\nu(d)$.

As before, let
${\cal Q}(\tau)\subset
{\cal Q}^1(S)$ be the set of all area one
quadratic differentials whose
horizontal measured geodesic lamination
is contained in ${\cal V}_0(\tau)$ and whose
vertical measured geodesic lamination hits $\tau$
efficiently.
As in Section 4, call $\tau\in {\cal T\cal T}$ \emph{$R-\xi$-tight}
for a number $R>0$ and some $\xi\in {\cal P\cal E}(\tau)$
if the diameter of
${\cal Q}(\tau)\cap W^{u}(\xi)$ with respect
to the lift $d$ of the Teichm\"uller metric
is at most $R$. For $\lambda\in {\cal C\cal L}(\tau)$
denote by $\{\tau_i(\lambda)\}$ the full splitting
sequence issuing from $\tau_0(\lambda)=\tau$ and determined by
$\lambda$. Let $m>0$ be the number of orbits of the action of
${\cal M}(S)$ on ${\cal T\cal T}$. The next technical observation
is the main ingredient for the
completion of the proof of the theorem from the
introduction. It roughly says that controlled recurrence implies
tightness.

\begin{lemma}\label{tightcon}
Let $\tau\in {\cal T\cal T}$ and let
$q_1\in{\cal Q}(\tau)$ be the lift of a forward recurrent
point $q_0\in {\cal Q}(S)$. Let $\xi$ be the vertical
projective measured geodesic lamination of $q_1$.
There are numbers
$R>2,k>4R+ma$ depending on $q_1$ and for every
complete geodesic lamination
$\lambda\in {\cal C\cal L}(\tau)\cap {\cal H}$
with $E(\lambda)\in {\cal R\cal M\cal L}(q_0)$ there
is a sequence $t_i\to \infty$ with the following
properties.
\begin{enumerate}
\item For each $i$ there are numbers $j(i)>0$, $\ell(i)>j(i)$
and there are numbers
$s\in [t_i,t_i+ma],t\in [t_i+k,t_i+k+ma]$
such that $e^s\nu(\lambda,\tau)\in
{\cal V}_0(\tau_{j(i)}(\lambda))$, $e^t\nu(\lambda,\tau)\in
{\cal V}_0(\tau_{\ell(i)}(\lambda))$
and such that the train tracks
$\tau_{j(i)}(\lambda),\tau_{\ell(i)}(\lambda)$
are $R-\xi$-tight.
\item For every $i$ there is some $g_i\in {\cal M}(S)$
such that $g_i\tau_{j(i)}(\lambda)=\tau_{\ell(i)}(\lambda)$.
\end{enumerate}
\end{lemma}
\begin{proof}
By Lemma \ref{simple2}
and inequality (\ref{Lipschitz}) in Section 2, for
all $q,z\in {\cal Q}^1(S)$ with $d(Pq,Pz)\leq a$
(where $a>0$ is as above), the distance in
${\cal C}(S)$ between $\Upsilon_{\cal Q}(q)$ and $\Upsilon_{\cal
Q}(z)$ is bounded from above by a universal constant $b>0$. Let
moreover $p>0$ be the maximal distance in ${\cal C}(S)$ between
any two vertex cycles of any complete train tracks $\eta_1,\eta_2$
on $S$ so that either $\eta_1=\eta_2$ or that $\eta_2$ can be
obtained from $\eta_1$ by a full split. Let $\chi_2 >0$ be as in
Lemma \ref{upsilon} and let $\ell >0$ be as in Lemma \ref{vertexcycle}.

By assumption, the $\omega$-limit set of the $\Phi^t$-orbit of
$q_0$ contains $q_0$. By Lemma 2.2 and the results of
\cite{H06a} we have
$d(\Upsilon_{\cal
Q}(\Phi^tq_1),\Upsilon_{\cal Q}(q_1))\to \infty$ $(t\to \infty)$.
Thus we can find small open relative compact
neighborhoods $V\subset U$ of $q_1$ in
${\cal Q}^1(S)$, numbers
$T_2>T_1>T_0=0$ and mapping classes $g_i\in {\cal M}(S)$ such that
for all $q\in V$ we have $P\Phi^{T_i}q\in g_iU$ and
\begin{equation}\label{ups} d(\Upsilon_{\cal Q}(\Phi^{T_{i}}q),
\Upsilon_{\cal Q}(\Phi^{T_{i-1}}q))\geq 2mp+
\ell+2\chi_2+2b \quad (i=1,2).
\end{equation}

Let for the moment $q\in V$ be an arbitrary
quadratic differential with $\pi(q)\in E({\cal H})\subset
{\cal F\cal M\cal L}$
and denote by $q_h,q_v$ the horizontal
measured geodesic lamination and the vertical
measured geodesic lamination of $q$, respectively.
Let $\sigma\in {\cal T\cal T}$ be a train track
which carries $q_h$ and hits $q_v$ efficiently.
We assume that there is some
$s_0\in [0,a]$ with $\Phi^{s_0}q\in {\cal Q}(\sigma)$.
Assume moreover that there is a complete geodesic lamination
$\lambda\in {\cal C\cal L}(\sigma)$ such that
$E(\lambda)$ is the projective class of $q_h$.
Then $\lambda$ determines a full splitting
sequence $\{\sigma_i(\lambda)\}$ issuing from
$\sigma_0(\lambda)=\sigma$.

For $i=1,2$ and for
$T_i>0$ as above
we can find some
$j_2>j_1>j_0=0,s_i\in [0,a]$ such that $e^{T_i+s_i}q_h\in
{\cal V}_0(\sigma_{j_i}(\lambda))$, which is equivalent
to saying that
$\Phi^{T_i+s_i}q\in {\cal Q}(\sigma_{j_i}(\lambda))$.
By our choice of $T_i$ we have
\begin{equation}
d(\Upsilon_{\cal Q}(\Phi^{T_i+s_i}q),
\Upsilon_{\cal Q}(\Phi^{T_{i-1}+s_{i-1}}q))\geq
2mp +\ell +2\chi_2 \quad(i=1,2).
\end{equation}
Thus by Lemma \ref{upsilon} and our choice of $p$,
for all $\alpha_i\in [j_i,j_i+m]$
the distance in ${\cal C}(S)$ between any
vertex cycle of $\sigma_{\alpha_{i}}(\lambda)$ and any
vertex cycle of $\sigma_{\alpha_{i-1}}(\lambda)$
is at least $\ell$.
Moreover, there is some $\tau_i\in [T_{i},T_{i}+ma]$ such that
$e^{\tau_i}q_h\in {\cal V}_0(\sigma_{\alpha_{i}}(\lambda))$ and
hence the total mass that the horizontal measured
geodesic lamination $e^{T_{i-1}}q_h$ disposes on the
complete train track $\sigma_{\alpha_{i}}(\lambda)$ is bounded
from below by $e^{-(T_{i}-T_{i-1})-ma}$.
Therefore Corollary \ref{tight},
applied to the train tracks $\sigma_{\alpha_2}(\lambda)\prec
\sigma_{\alpha_1}(\lambda)\prec\sigma$,
shows the existence
of a universal constant $R>0$ such that the
train track $\sigma_{\alpha_1}(\lambda)$
is $R-[q_v]$-tight where $[q_v]$ is the projective
class of the vertical measured geodesic lamination of $q$.
Note that the number $R>0$ only depends on the point $q_1$
but not on $q\in V$ or
the train track $\sigma\in {\cal T\cal T}$ with
$\Phi^{s_0}q\in {\cal Q}(\sigma)$.

For this number $R>0$, choose a
number $T_3>T_2+4R+2ma$ such that
\begin{equation}
d(\Upsilon_{\cal Q}(\Phi^{T_3}q_1),\Upsilon_{\cal Q}(\Phi^{T_2}q_1))
\geq 2mp+\ell +2\chi_2+b\end{equation}
and that $\Phi^{T_3}q_1\in g_3V$ for some
$g_3\in {\cal M}(S)$; such a
number exists since the projection
$q_0$ of $q_1$ to ${\cal Q}(S)$
is forward recurrent. Choose an open neighborhood
$W\subset V$ of $q_1$ such that $\Phi^{T_i}q\in g_iV$ for all
$q\in W$ and $i=1,2,3$.
Write $T_4=T_3+T_1$ and $T_5=T_3+T_2$; then we have
$\Phi^{T_j}q\in\cup_{h\in {\cal M}(S)}hU$ for all $q\in W$ 
and every $j\in \{0,\dots,5\}$.
Define $k=T_4-T_1\geq 4R+2ma$ 
and note that $k$ only depends on $q_1$.

As above, let $\tau\in {\cal T\cal T}$ be a train track which
carries the horizontal projective measured geodesic lamination of
$q_1$ and such that the vertical projective measured geodesic
lamination $\xi$ of $q_1$ hits $\tau$ efficiently. Let $\lambda\in
{\cal C\cal L}(\tau)\cap {\cal H}$
be such that $E(\lambda)\in {\cal R\cal
M\cal L}(q_0)$. Let $q\in {\cal Q}(\tau)\cap W^u(\xi)$ be such
that the horizontal projective measured geodesic lamination of $q$
equals $E(\lambda)$. For the neighborhood $W\subset {\cal Q}^1(S)$
of $q_1$ as above, the set $\{t>0\mid \Phi^tq\in \tilde
W=\cup_{h\in {\cal M}(S)} hW\}$ is unbounded. Let
$\{r_n\}\subset [0,\infty)$ be a sequence tending to infinity
such that $\Phi^{r_n}q\in h_nW$ for some $h_n\in {\cal M}(S)$
and every $n>0$.

Let $\{\tau_i(\lambda)\}$ be the full splitting sequence determined by
$\tau=\tau_0(\lambda)$ and $\lambda$. For $n>0$
we have $\Phi^{r_n}q\in h_nW$.
There is a number $s_0\in [0,a]$ and a number
$j_0(n)>0$ such that $\Phi^{t_n+s_0}q\in {\cal Q}(\tau_{j_0(n)}(\lambda))$.
Using the above constants $T_i>0$, there are numbers
$j_5(n)>j_4(n)>j_3(n)>j_2(n)>j_1(n)>j_0(n)$ 
and numbers $s_i\in [0,a]$ such that
$\Phi^{T_i+s_i+r_n}q
\in {\cal Q}(\tau_{j_i(n)}(\lambda))$ $(i=1,2,3,4,5)$.
Since there are only $m$ distinct orbits of complete train tracks
under the action of the mapping class group, there
are moreover numbers $j(n)\in [j_1(n),j_1(n)+m]$
and $\ell(n)\in [j_4(n),j_4(n)+m]$ and
there is some
$g\in {\cal M}(S)$ such that
$g\tau_{j(n)}(\lambda)=\tau_{\ell(n)}(\lambda)$. By the above consideration,
the train tracks $\tau_{j(n)}(\lambda), \tau_{\ell(n)}(\lambda)$ are
$R-\xi$-tight. Then the sequence $t_n=r_n+T_1$ and the numbers 
$j(n)>0,\ell(n)>j(n)$ have the
required properties stated in the lemma with $k=T_4-T_1\geq 4R+ma$.
\end{proof}

Now we are ready to complete the main step in
the proof of the theorem from the introduction.

\begin{proposition}
\label{fillingisleb}  An ${\cal M}(S)$-invariant
Radon measure on ${\cal M\cal L}$ which
gives full mass to the filling measured geodesic laminations
coincides with the Lebesgue measure up to scale.
\end{proposition}
\begin{proof}
By Lemma \ref{groupnontrivial} and
Lemma \ref{fillingrec}, we only have to show
that there is no ${\cal M}(S)$-invariant ergodic Radon
measure on ${\cal M\cal L}$ with $H_\eta=\{0\}$ which
gives full mass to the recurrent measured geodesic laminations.

For this we argue by contradiction and we assume that
such a Radon measure $\eta$ exists. Using once
more the cocycle reduction theorem of Sarig,
there is a Borel function $u:{\cal P\cal M\cal L}\to \mathbb{R}$
such that the measure $\eta$
gives full mass to the graph $\{(x,u(x))\mid
x\in {\cal P\cal M\cal L}\}$ of the function $u$.
In particular, $\eta$ projects to an ${\cal M}(S)$-invariant
measure class $\hat \eta$ on ${\cal P\cal M\cal L}$ which
gives full mass to the set ${\cal R\cal M\cal L}$ of
recurrent points.
Using the notations from Proposition \ref{Vitali} and its
proof, this means that there is a forward recurrent point
$q_0\in {\cal Q}(S)$ such that the measure class
$\hat \eta$ on ${\cal P\cal M\cal L}$ gives full measure to the set
${\cal R\cal M\cal L}(q_0)$.

To derive a contradiction we adapt the arguments
of Ledrappier and Sarig \cite{LS06} to our
situation. Denote again by ${\cal H}\subset {\cal C\cal L}$
the set of all complete geodesic laminations
containig a minimal component which fill up $S$ and let
$E:{\cal H}\to {\cal F\cal M\cal L}$ be as before.
Every finite Borel measure
$\mu$ on ${\cal P\cal M\cal L}$ which gives full
mass to the returning projective measured geodesic laminations
induces a finite Borel measure $\tilde \mu$
on ${\cal H}$.
Namely, returning projective measured geodesic laminations
are uniquely ergodic and hence
for a Borel subset $C$ of ${\cal H}$ we can define
\begin{equation}
\tilde \mu(C)=\int_{E(C)}
\sharp (E^{-1}(z)\cap C)d\mu(z).
\end{equation}

Let again $q_0\in {\cal Q}(S)$ be a forward recurrent
quadratic differential such that
the measure class $\hat \eta$ gives full mass
to ${\cal R\cal M\cal L}(q_0)$. Let $q_1\in {\cal Q}^1(S)$
be a lift of $q_0$. We may assume that the
vertical measured geodesic lamination $\xi=\pi(-q_1)$
of $q_1$ is uniquely ergodic and fills up $S$.
We may moreover assume that there is a train track
$\tau\in {\cal T\cal T}$ with $q_1\in {\cal
Q}(\tau)$ (as before, this can
be achieved by possibly replacing $q_1$ by $\Phi^tq_1$ for
some $t>0$). By the considerations in Section 4 we may 
assume that $\tau$ is $R_0-\xi$-tight for a number
$R_0>0$.

Let $\rho\in
\mathbb{R}$ be such that $\cup_{\rho-1\leq s\leq \rho+1}\Phi^s
W^{su}(q_1)$ contains some density point of
the locally finite measure $\eta^u$ on $W^u(q_1)$ induced from $\eta$
whose
horizontal measured geodesic lamination is carried by the train
track $\tau$; such a number exists since by invariance under the
mapping class group, the ${\cal M}(S)$-invariant measure class
$\hat\eta$ on ${\cal P\cal M\cal L}$ is of full support and since
moreover the set of all measured geodesic laminations
carried by $\tau$ has non-empty interior.

Let $\hat \mu_0$ be the restriction of the measure
$\eta^u$ to the
set $\cup_{\rho-1\leq s\leq \rho+1} \Phi^s W^{su}(q_1)\subset W^u(q_1)$
and denote by $\mu_0$ the projection of $\hat \mu_0$ to ${\cal
P\cal M\cal L}$. Since $\eta$ is a Radon measure on ${\cal M\cal
L}$ by assumption, the measure $\mu_0$ is a locally finite Borel
measure on $\pi W^u(\xi)\subset {\cal P\cal M\cal L}-\xi$ which
gives full mass to ${\cal R\cal M\cal L}(q_0)$.
Since $\tau$ is $R_0-\xi$-tight for some $R_0>0$,
the set ${\cal Q}(\tau)\cap W^{u}(\xi)$ is relative compact and therefore
the intersection of the support of $\hat \mu_0$ with
$\cup_t\Phi^t{\cal Q}(\tau)$ is compact as well.
Since $\eta$ and hence $\hat \mu_0$ is Radon,  
the total $\mu_0$-mass of the set of all projective
measured geodesic laminations which are carried by the
train track $\tau$ is finite.
By the above consideration, $\mu_0$ induces a 
finite nontrivial Borel measure $\tilde \mu_0$ on ${\cal
C\cal L}(\tau)$ which gives full mass to the set of complete
geodesic laminations $\zeta\in {\cal H}\cap
{\cal C\cal L}(\tau)$ with $E(\zeta)\in {\cal
R\cal M\cal L}(q_0)$. 

Similarly, for the constants $R>2,m>0,k>4R+ma$ as in
Lemma \ref{tightcon},
let $\tilde \mu_1$ be the finite Borel measure on
${\cal C\cal L}(\tau)$ which is induced from the restriction of
$\eta^u$ to $\cup_{\rho+2R\leq s\leq \rho+k+2R+ma} \Phi^sW^{su}(q_1)$.
Since $H_{\eta}=\{0\}$ by assumption, the measures $\tilde
\mu_0,\tilde \mu_1$ are singular.

Define a \emph{cylinder} in ${\cal C\cal L}(\tau)$
to be a set of the
form ${\cal C\cal L}(\sigma)$ where
$\sigma\in {\cal T\cal T}$ is a complete train track
which can be obtained from $\tau$ by a full
splitting sequence. A cylinder is a subset
of ${\cal C\cal L}$ which is both open and closed
\cite{H06b}. The intersection of two cylinders
is again a cylinder. Since every
point in ${\cal C\cal L}(\tau)$
is an intersection of countably many cylinders,
the $\sigma$-algebra on ${\cal C\cal L}(\tau)$
generated by cylinders is
the usual Borel $\sigma$-algebra.
Now $\tilde \mu_0,\tilde \mu_1$ are mutually singular
Borel measures on ${\cal C\cal L}(\tau)$ and hence there is
a cylinder ${\cal C\cal L}(\sigma)\subset {\cal C\cal L}(\tau)$
such that $\tilde \mu_0({\cal C\cal L}(\sigma))>
2\tilde \mu_1({\cal C\cal L}(\sigma))$.

By Lemma \ref{tightcon}, for $\tilde \mu_0$-almost
every $\lambda\in {\cal C\cal L}(\sigma)$ there is
a sequence $j(i)\to \infty$ and there is some
$g(i,\lambda)\in {\cal M}(S)$ with the following properties.
\begin{enumerate}
\item The train tracks
$\tau_{j(i)}(\lambda),g(i,\lambda)\tau_{j(i)}(\lambda)$ are
$R-\xi$-tight.
\item $g(i,\lambda)
{\cal C\cal L}(\tau_{j(i)}(\lambda))\subset
{\cal C\cal L}(\tau_{j(i)}(\lambda))$.
\item For $q\in {\cal Q}(g(i,\lambda)\tau_{j(i)}(\lambda))\cap W^u(\xi)$
there is some $t\in [k-ma-2R,k+ma+2R]$ such that $e^{-t}q\in
{\cal Q}(\tau_{j(i)}(\lambda)) \cap W^u(\xi)$.
\end{enumerate}
These properties imply the following. Let $q\in \cup_{\rho-1\leq
t\leq \rho+1}\Phi^tW^{su}(q_1)$ be such that the horizontal measured
geodesic lamination $q_h$ of $q$ is carried by
$\tau_{j(i)}(\lambda)$. Let $s\in \mathbb{R}$ be such that
$\Phi^sq\in {\cal Q}(\tau_{j(i)}(\lambda))$ and let
$z=W^{ss}(g(i,\lambda)\Phi^sq)\cap W^u(\xi)$; then the horizontal measured
geodesic lamination of $z$ is carried by $\tau_{j(i)}(\lambda)$,
and we have $z\in \Phi^{t+s}W^{su}(q_1)$ for some $t\in
[k-ma-2R,k+ma+2R]$. Since $k>4R+ma$
it is now immediate
from invariance of the measure $\eta$ under the
action of ${\cal M}(S)$, from
the definitions of the measures
$\tilde \mu_0,\tilde \mu_1$ and the fact that
the action of ${\cal M}(S)$ on ${\cal M\cal L}$ commutes
with the action of the group of translations
that $\tilde \mu_0({\cal C\cal
L}(\tau_{j(i)}\lambda))\leq \tilde \mu_1({\cal C\cal
L}(\tau_{j(i)}\lambda))$.

On the other hand, there is a countable partition of a subset of
${\cal C\cal L}(\sigma)$ of full $\tilde \mu_0$-mass into
cylinders ${\cal C\cal L}(\sigma_i)$ with train tracks
$\sigma_i\in {\cal T\cal T}$ $(i>0)$ which can be obtained from
$\sigma$ by a full splitting sequence and which satisfy 1),2),3)
above. This partition can inductively be constructed as follows.
Beginning with the train track $\sigma$, there is a full splitting
sequence of minimal length $n\geq 0$ issuing from $\sigma$ which
connects $\sigma$ to some train track $\sigma_1\in {\cal T\cal T}$
with the above properties. Let $\eta_1,\dots,\eta_k$ be the
collection of all train tracks which can be obtained from $\sigma$
by a full splitting sequence of length $n$ and assume after
reordering that we have $\sigma_1=\eta_1$. Repeat this construction
simultaneously with the train tracks $\eta_2,\dots,\eta_k$. After
countably many steps we obtain a partition of $\tilde
\mu_0$-almost all of ${\cal C\cal L}(\sigma)$ as required.

Together we conclude that necessarily $\tilde \mu_1({\cal C\cal
L}(\sigma))\geq \tilde \mu_0({\cal C\cal L}(\sigma))$ which
contradicts our choice of $\sigma$. In other words,
the case $H_{\eta}=\{0\}$
is impossible which completes the proof of our proposition.
\end{proof}

{\bf Remark:} Let $q_0\in {\cal Q}(S)$ be a forward recurrent
point.
The arguments in the proof of Proposition \ref{fillingisleb}
can be used to construct a Vitali relation
for the lift to ${\cal C\cal L}$
of any Borel measure on ${\cal R\cal M\cal L}(q_0)$.
In other words, with some extra arguments, Proposition
\ref{Vitali} can be deduced from
Proposition \ref{fillingisleb} and its proof.
However, we
included Proposition \ref{Vitali} in the present form
since its basic idea is simpler and more
geometric, moreover it is used
in \cite{H07}.

Choose a complete hyperbolic metric on $S$ of finite volume.
Let $S_0$ be a proper connected bordered subsurface of $S$
with geodesic boundary. Then $S_0$ has
negative Euler characteristic.
We allow that distinct
boundary components of $S_0$ are defined by the
same simple closed geodesic in $S$.
Denote by ${\cal M}(S_0)$ (or ${\cal M}(S-S_0)$)
the subgroup of ${\cal M}(S)$ of all
elements which can be represented
by a diffeomorphism fixing $S-S_0$ (or $S_0$) pointwise.
The stabilizer
${\rm Stab}(S_0)$ of $S_0$
in ${\cal M}(S)$ contains a subgroup of finite index of
the form
${\cal M}(S_0)\times {\cal M}(S-S_0)\times
{\cal D}(\partial S_0)$ where
${\cal D}(\partial S_0)$ is the free abelian group of
Dehn twists about the geodesics in $S$
which define the boundary of $S_0$.

Let $\widehat S_0$
(or $\widehat{S-S_0}$) be the
surface of finite type which we obtain from
$S_0$ (or $S-S_0)$ by
collapsing each boundary circle to a puncture.
The space ${\cal M\cal L}(S_0)$ of all
measured geodesic laminations on $\widehat S_0$
can be identified with the space of all measured
geodesic laminations on $S$ whose support
is contained in the interior of $S_0$.
We say that a measured geodesic
lamination $\nu\in {\cal M\cal L}(S_0)$ \emph{fills}
$S_0$ if its support is minimal and intersects
every simple closed geodesic contained in the
interior of $S_0$ transversely.
There is a ray of
Lebesgue measures on ${\cal M\cal L}(S_0)$ which
are invariant under the mapping class group
${\cal M}(\widehat S_0)>{\cal M}(S_0)$ of the
surface $\widehat S_0$ and hence it is invariant
under ${\rm Stab}(S_0)$.
For every $\phi\in {\cal M}(S)-{\rm Stab}(S_0)$,
the image $\phi(S_0)$ of $S_0$ under $\phi$ is a subsurface
of $S$ which is distinct from $S_0$.
The image $\phi(\zeta)$ under $\phi$
of a measured geodesic lamination
$\zeta$ on $S_0$ which fills $S_0$
is a measured geodesic lamination
which fills $\phi(S_0)$ and hence
this image is not contained in ${\cal M\cal L}(S_0)$.
Note that we have ${\rm Stab}(\phi(S_0))=
\phi\circ {\rm Stab}(S_0)\circ \phi^{-1}$.

Now let $c$ be any (possibly trivial)
simple weighted geodesic multicurve
on $S$ which is disjoint from the interior of $S_0$.
Then for every $\zeta\in {\cal M\cal L}(S_0)$ the
union $c\cup \zeta$ is a measured geodesic lamination
on $S$ in a natural way which we denote by
$c\times \zeta$. Thus $c\times {\cal M\cal L}(S_0)$
is naturally a closed subspace of ${\cal M\cal L}$. This
subspace
can be equipped with a
${\rm Stab}(c\cup S_0)< {\rm Stab}(S_0)$-invariant
ergodic Radon measure
$\mu_{c,S_0}$ induced by a measure $\mu_0$ from our
ray of ${\cal M}(\widehat S_0)$-invariant
Lebesgue measures on ${\cal M\cal L}(S_0)$.
By invariance, we obtain
a ${\cal M}(S)$-invariant ergodic wandering measure
on ${\cal M\cal L}$ by defining
\begin{equation}
\lambda_{c\times S_0}=\sum_{\phi\in {\cal M}(S)}\phi_*\mu_{c,S_0}.
\end{equation}
We call $\lambda_{c\times S_0}$ a \emph{standard subsurface
measure} of $S_0$.
If the support of the weighted multi-curve $c$ contains every
boundary component
of $S_0$ then we call the resulting ${\cal M}(S)$-invariant
measure $\lambda_{c\times S_0}$ on ${\cal M\cal L}$ a
\emph{special standard subsurface measure} on
${\cal M\cal L}$.

Recall from the introduction that a
\emph{rational} ${\cal M}(S)$-invariant measure
on ${\cal M\cal L}$ is a sum of weigthed Dirac masses
supported on the orbit of a simple weighted multi-curve.
Such a rational measure is a special standard subsurface
measure on ${\cal M\cal L}$ (for the empty subsurface).
We have.

\begin{lemma}\label{localfinite}
A special standard subsurface measure on ${\cal M\cal L}$
is locally finite.
\end{lemma}
\begin{proof}
Let $g$ be any complete hyperbolic metric on $S$
of finite volume. Then for every measured
geodesic lamination $\mu$ on $S$ the
\emph{$g$-length} $\ell_g(\mu)$
of $\mu$ is defined. By definition,
this length is the total mass of the measure on
$S$ which is the product of the transverse measure
for $\mu$ and the hyperbolic length element on
the geodesics contained in the support of $\mu$. For every
\emph{compact} subset $K$ of ${\cal M\cal L}$ there
is a number $m>0$ such that $K\subset K(m)=
\{\mu\in{\cal M\cal L}\mid \ell_g(\mu)\leq m\}$.

To show the lemma observe first that
for every weighted geodesic multi-curve $c$ on $S$ and
every $m>0$ the
set $K(m)$ contains only finitely many images of $c$
under the action of the mapping class group.
Namely, let $a>0$ be the minimal weight of a component
of the support of $c$. Then for every $\phi\in {\cal M}(S)$ the
$g$-length of the multi-curve $\phi(c)$ is not smaller than
$a$ times the maximal length of any closed geodesic
on the hyperbolic surface
$(S,g)$ which is freely homotopic to a component of $\phi(c)$.
However, there are only finitely many simple
closed geodesics on $S$ whose
$g$-length is at most $m/a$ and hence
the intersection of $K(m)$ with the ${\cal M}(S)$-orbit
of $c$ is indeed finite. In particular, a rational
${\cal M}(S)$-invariant measure on ${\cal M\cal L}$
is Radon.

Now let $S_0$ be a bordered connected subsurface of $S$
of negative Euler characteristic and
geodesic boundary. Let $c$ be a simple weighted
geodesic multi-curve which contains every boundary
component of $S_0$. Then the stabilizer ${\rm Stab}(c)$
of $c$ in ${\cal M}(S)$ contains the stabilizer
${\rm Stab}(c\cup S_0)$ of $c\cup S_0$ as a subgroup of
finite index. In particular, a
${\cal M}(\widehat{S_0})$-invariant
measure on ${\cal M\cal L}(S_0)$ in the Lebesgue
measure class induces an ergodic
${\rm Stab}(c)$-invariant
Radon measure on the space of laminations containing
$c$ as a component.

On the other hand, if $\phi \in {\cal M}(S)$
\emph{does not} stabilize
$c$ then $\phi$
moves at least one component of $c$ away from $c$.
Therefore by the above consideration, there are
only \emph{finitely many} cosets in
${\cal M}(S)/{\rm Stab}(c)$
containing some representative
$\phi$ such that
$\phi(c\times {\cal M\cal L}(S_0))\cap K(m)\not=\emptyset$.
This shows that a special subsurface measure on
${\cal M\cal L}$ is Radon.
\end{proof}

Finally we are able to complete the proof of
the theorem from the introduction.

\begin{theorem}
Let $\eta$ be an ${\cal M}(S)$-invariant
ergodic Radon
measure on the space ${\cal M\cal L}$
of all measured geodesic laminations on $S$.
\begin{enumerate}
\item If $\eta$ is non-wandering then $\eta$
is the Lebesgue measure up to
scale.
\item If $\eta$ is wandering then
either $\eta$ is rational or $\eta$
is a standard subsurface measure.
\end{enumerate}
\end{theorem}
\begin{proof}
Let $\eta$ be an ergodic ${\cal M}(S)$-invariant
Radon measure on ${\cal M\cal L}$. By Proposition
\ref{fillingisleb}, if $\eta$ gives full mass to the
measured geodesic laminations which fill up $S$
then $\eta$ coincides with the Lebesgue measure
up to scale. Thus by ergodicity
and invariance we may assume that $\eta$ gives
full mass to the measured geodesic laminations
which do not fill up $S$.

Let $\lambda$ be a density point for $\eta$.
The support of $\lambda$ is a union of components
$\lambda_1\cup\dots\cup \lambda_k$. We
assume that
these components are ordered in such a way that
there is some $\ell\leq k$ such that the
components $\lambda_1,\dots,\lambda_{\ell}$ are
minimal arational and that the components
$\lambda_{\ell+1},\dots,\lambda_k$ are simple closed curves.
By ergodicity, the same decomposition then holds for
$\eta$-almost every $\lambda\in{\cal M\cal L}$.

If $\ell=0$ then by ergodicity, $\eta$ is rational and there
is nothing to show. Thus assume that
$\ell>0$. Then $\mu=\lambda_1\cup \dots\cup \lambda_\ell$
\emph{fills} a subsurface $S_0$ of $S$ with $\ell$
connected components, each of which is of negative
Euler characteristic. In particular, by ergodicity the
measure $\eta$ is wandering and gives
full mass to the set $\{g(c\times {\cal M\cal L}(S_0))\mid
g\in {\cal M}(S)\}$.

Write $c=\lambda_{\ell +1}\cup \dots\cup \lambda_k$.
Since $\lambda$ is a density
point for $\eta$, the restriction of $\eta$ to
$c\times {\cal M\cal L}(S_0)$ does not  vanish.
However, this restriction is a ${\rm Stab}(c\cup S_0)$-invariant
Radon measure on ${\cal M}(S_0)$ which gives full mass to
the measured geodesic laminations filling up $S_0$ and
therefore this restriction is an interior point of
our cone of ${\rm Stab}(c\cup S_0)$-invariant
Lebesgue measure on $c\times {\cal M\cal L}(S_0)$.
In other words, $\eta$ is a standard
subsurface measure. This completes the proof of the
theorem.
\end{proof}

{\bf Remark:} In \cite{LM07}, Lindenstrauss and
Mirzakhani obtain a stronger result.
They show that a locally finite
standard
subsurface measure on ${\cal M\cal L}$ is special.

\section*{Appendix}

The purpose of this appendix is to present some
results from
the paper \cite{MW02} of Minsky and Weiss
in the form needed in Section 3.

As in Section 2, denote by ${\cal Q}(S)$ the moduli
space of area one
holomorphic quadratic differentials on $S$.
Every $q\in {\cal Q}(S)$ defines an isometry
class of a singular
euclidean metric on $S$. The set $\Sigma$
of singular points for this metric
coincides precisely with the set of zeros for
$q$. We also assume that the differential
has a simple pole at each of the punctures of $S$ and
hence it can be viewed as a meromorphic
quadratic differential
on the compactified surface $\hat S$ which we obtain
by filling in the punctures in the standard way.

A \emph{saddle connection} for $q$ is a path
$\delta:(0,1)\to S-\Sigma$ whose image in each
chart is an euclidean straight line and which
extends continuously to a path
$\bar \delta:[0,1]\to \hat S$ mapping the endpoints to
singularities or punctures.
A saddle connection
does not have self-intersections.
Two saddle connections
$\delta_1,\delta_2$ are \emph{disjoint} if
$\delta_1(0,1)\cap \delta_2(0,1)=\emptyset.$
The closure of any finite collection of
pairwise disjoing saddle connections
on $S$ is an embedded graph in $S$.
By Proposition 4.7 of \cite{MW02} (see also \cite{KMS86}),
the number of pairwise disjoint saddle connections
for a quadratic differential $q\in {\cal Q}(S)$ is
bounded from above by a universal constant
$M>0$ only depending on the topology of $S$.

Recall that a \emph{tree} is a graph without
circuits.
For $\epsilon >0$ let $K(\epsilon)\subset {\cal Q}(S)$
be the set of all quadratic differentials $q$ such that
the collection of all saddle connections of
$q$ of length at most $\epsilon$ is a tree.
We have.

\begin{lem}
For every $\epsilon >0$ the set
$K(\epsilon)\subset {\cal Q}(S)$ is compact.
\end{lem}
\begin{proof}
It is enough to show that for every $q\in K(\epsilon)$
the $q$-length of any
simple closed curve on $S$ is bounded from below
by $\epsilon$ (see \cite{R05,R06}).

Thus let $c$ be any simple closed geodesic on
$S$ for the $q$-metric. Then up to replacing
$c$ by a freely homotopic simple closed curve of the
same length we may assume that
$c$ consists of a sequence of saddle connections
for $q$. Since the set of saddle connections
of length at most $\epsilon$ does not contain
a circuit, the curve $c$ contains at least one
saddle connection of length at least $\epsilon$.
But this just means that the $q$-length of $c$
is at least $\epsilon$ as claimed.
\end{proof}

The following proposition is a modified version
of Theorem 6.3 of \cite{MW02}. We use the
notations from \cite{MW02}. Let
${\cal L}_q$ be the set of all saddle connections
of the quadratic differential $q$.
For $k\geq 1$ define
\[{\cal E}_k=\{E\subset {\cal L}_q\mid E\, \text{consists of $k$ disjoint
segments}\}.\]
Denote again
by $h_t$ the horocycle flow on ${\cal Q}(S)$.
For $E\in {\cal E}_k$ and $t\in \mathbb{R}$ define
$\ell_{q,E}(t)=\max_{\delta\in E}\ell_{q,\delta}(t)$ where
$\ell_{q,\delta}(t)$ is the length of $\delta$ with respect to the
singular euclidean metric defined by $h_tq$. For $k\geq 0$ let
\[\alpha_{k}(t)=\min_{E\in {\cal E}_{k}}\ell_{q,E}(t).\]

\begin{prop}
There are positive constants $C,\alpha,\rho_0$ depending
only on $S$ with the following property.
Let $q\in {\cal Q}(S)$, let $I\subset \mathbb{R}$
be an interval and let $0<\rho^\prime\leq\rho_0$.
Define \[A=\{\delta\in {\cal L}_q\mid
\ell_{q,\delta}(t)\leq \rho^\prime\}\quad\text{for all } t\in I.\]
If $\cup\{\overline{\delta}\mid \delta \in A\}\subset S$ is an
embedded tree with $r\geq 0$ edges then
for any $0<\epsilon <\rho^\prime$ we have:
\[\vert \{t\in I\mid \alpha_{r+1}(t)<\epsilon\}\vert \leq
C\bigl(\frac{\epsilon}{\rho^\prime}\bigl)^\alpha\vert I\vert.\]
\end{prop}
\begin{proof}
Let $M>0$ be such that for every $q\in {\cal Q}(S)$
the number of pairwise disjoint
saddle connections of $q$
is bounded from above by $M-1$. By Proposition 6.1 of \cite{MW02}
there is a number $\rho_0>0$ with the following
property. If $E\in {\cal E}_k$ is such that the closure $S(E)$ of the
union of all simply connected components of
$S-\cup_{\delta\in E}\bar\delta$ is all of $S$ then
$\ell_{q,E}(0)\geq\rho_0$.

Let $q\in {\cal Q}(S)$ and let $A\subset {\cal L}_q$ be a union
of pairwise disjoint saddle connections whose closure is
an embedded graph in $S$ without circuits.
Assume that $A$ consists of $r\geq 0$ segments. We necessarily
have $r<M$. For a number $C>0$ to be determined later
let $0<\epsilon <C\rho^\prime$ and let
\[V_\epsilon =\{t\in I\mid \alpha_{r+1}(t)<\epsilon\}.\]
For $k=1,\dots,M-r-1$ define
\[L_k=\epsilon \bigl(\frac{\rho^\prime}{\epsilon}
\bigr)^{\frac{k-1}{M-r-1}}.\]
We choose $C>0$ in such a way that
$L_k/L_{k+1}\leq C^{\frac{1}{M-r-1}}$.
For $t\in V_\epsilon$ let
\[\kappa(t)=\max\{k\mid \alpha_k(t)<L_k\}\]
and let $V_k=\{t\in V_\epsilon\mid \kappa(t)=k\}$.
Then $\kappa(t)\leq M-1$
for all $t$ and hence
$V_\epsilon$ is the disjoint union of the measurable
sets $V_{k+1},\dots, V_{M-1}$ .
Thus there
is some $k\in \{r+1,\dots, M-1\}$ for which
\[\vert V_k\vert \geq \frac{\vert V_\epsilon \vert}{M-r-1}.\]
For this choice of $k$ define $L=L_k$ and $U=L_{k+1}$.
Note that we have
\[\alpha_{\kappa(t)}(t)<L_{\kappa(t)},\quad
\alpha_{\kappa(t)+1}\geq L_{\kappa(t)+1}.\]

Following \cite{MW02}, for $\delta\in {\cal L}_\delta-A$
let $H(\delta)$ be the set of $t\in I$ for which
$\ell_{q,\delta}(t)<L$, and whenever
$\delta\cap \delta^\prime\not=\emptyset$ for
$\delta\not=\delta^\prime\in {\cal L}_q$ we have
\[\ell_{q,\delta^\prime}(t)\geq \frac{U\sqrt{2}}{3}.\]

Following the argument in Section 6 of \cite{MW02} we only have to
verify that $V_k\subset \cup_{\delta\in {\cal L}_q-A}
H(\delta)$. Namely, let $t\in V_k$ and let $E\in {\cal E}_k$
be such that $\ell_{q,E}(t)=\alpha_k(t)<L$.
Denote by $S(E)$ the closure of the union
of the simply connected components of $S-\cup_{\delta\in E}\bar\delta$.
By Proposition 6.1 of \cite{MW02} we have
$S(E)\not=E$ and hence since $k>r$ and the graph
defined by the saddle connections contained in $A$ does not have
circuits, the boundary of $S(E)$ contains
at least one saddle connection $\delta$
which is \emph{not} contained in $A$.
But this just means that $t\in H(\delta)$ (see Claim 6.7 in \cite{MW02}).
This complete the proof of the proposition.
\end{proof}

As in \cite{MW02} we use the lemma and the
proposition to derive
a recurrency property for the horocycle flow.
For its formulation, denote by $\chi_C$ the characteristic function of
the set $C\subset {\cal Q}(S)$.

\begin{theo}
For any $\epsilon >0$ there is a compact
set $K\subset {\cal Q}(S)$ such that for
any $q\in {\cal Q}(S)$ with minimal horizontal
measured geodesic lamination which fills $S$
we have
\[{\rm Avg}_{t,q}(K)=\lim\inf_{t\to \infty}
\frac{1}{t}\int_0^t \chi_K(h_tq)dt \geq 1-\epsilon.\]
\end{theo}
\begin{proof}
Let $q$ be a quadratic differential with horizontal
measured geodesic lamination which fills up $S$.
Then the horizontal saddle connections of $q$ form
an embedded graph without circuits. Morever, the
number of these saddle connections is bounded
from above by a universal constant. Now if $\delta$
is any saddle connection whose length is constant
along the horocycle flow then $\delta$ is horizontal.
But this just
means that we can apply the above proposition
as in the proof of Theorem H2 of \cite{MW02} to
obtain the theorem.
\end{proof}

\bigskip

\noindent
MATHEMATISCHES INSTITUT DER UNIVERSIT\"AT BONN,\\ 
BERINGSTRASSE 1,\\
D-53115 BONN, GERMANY\\

\smallskip
\noindent
e-mail: ursula@math.uni-bonn.de

\end{document}